\documentclass{article}
\usepackage{amsmath,amssymb,enumerate}
\usepackage{verbatim}
\usepackage[all]{xy}

\newtheorem{prop}{Proposition}[section]
\newtheorem{coro}[prop]{Corollary}
\newtheorem{lem}[prop]{Lemma}
\newtheorem{rem}[prop]{Remark}
\newtheorem{exe}[prop]{Example}
\newtheorem{defi}[prop]{Definition}
\newtheorem{theo}[prop]{Theorem}

\newtheorem{theoi}{Theorem}

\newenvironment{pf}{{\bfseries Proof:}}{\ \hfill $\Box$}

\def\D{\mathbb{D}}
\def\C{\mathbb{C}}
\def\R{\mathbb{R}}
\def\Q{\mathbb{Q}}
\def\Z{\mathbb{Z}}
\def\N{\mathbb{N}}

\def\V{\mathbb{V}}
\def\H{\mathbb{H}}

\def\M{\mathbb{M}}
\newcommand{\Qbar}{\overline{\Q}}

\newcommand{\Or}{\mathcal O_{\rho}}

\newcommand{\Ot}{\mathcal O_{T} }
\newcommand{\rep}{R(\pi_1(X,x), GL_N)}
\newcommand{\m}{\mathfrak m}

\newcommand{\hdot}{\bullet}

\begin{document}
\title{Linear Shafarevich Conjecture}
\author{P. Eyssidieux, L. Katzarkov\footnote{Partially supported by NSF Grant
    DMS0600800, by NSF FRG DMS-0652633, FWF grant P20778 and ERC grant.},
T. Pantev\footnote{Partially supported by NSF grant
DMS-0700446 and NSF Research Training Group Grant DMS-0636606.}, M. Ramachandran }

\date{April, 2, 2009}

\maketitle

\tableofcontents

\section*{Introduction}
\addcontentsline{toc}{section}{Introduction}

A complex analytic space $S$ is {\em holomorphically convex} if there
is a proper holomorphic morphism $\pi : S\to T$ with
$\pi_*\mathcal{O}_S=\mathcal{O}_T$ such that $T$ is a Stein space. $T$
is then called the {\em Cartan-Remmert reduction} of $S$.

The so-called Shafarevich conjecture of holomorphic convexity predicts
that the universal covering space $\widetilde{X^{\text{univ}}}$ of a
complex compact projective manifold $X$ should be holomorphically
convex.

This is trivial if the fundamental group is finite. The Shafarevich conjecture is a corollary of the
Riemann uniformization theorem in dimension $1$.  
The study   of  the  Shafarevich conjecture for
 smooth projective surfaces was initiated in the mid 80s by Gurjar and Shastri \cite{GUR} and  Napier
\cite{NAP}.

In the mid nineties, new  ideas introduced   by J.  K\'ollar and independently by F. Campana have
revolutionized the subject. The outcome was what is still the best result available 
with no assumption on the fundamental group, namely the 
construction of the Shafarevich map (aka $\Gamma$-reduction)
\cite{Cam1, Kol1, Kol2}.

At the same time, Corlette and Simpson \cite{Cor1, Cor2, Sim1, Sim2, Sim3}  were
developing Nonabelian Hodge theory. A bit later a  $p$- adic version of
Nonabelian Hodge theory in degree 1 was developed by Gromov and Schoen
\cite{GrS}.

The idea that Nonabelian Hodge theory can be used to
prove
Shafarevich conjecture was introduced in 1994 by the second author. He
 proved Shafarevich conjecture  for nilpotent
fundamental groups \cite{Katnilp}. At about the same time the second and the fourth author proved the Shafarevich conjecture for smooth projective
surfaces with fundamental group admitting  faithful Zariski dense representation,
 in a reductive complex algebraic group \cite{KR}.

The first author then found a way to extend this Nonabelian Hodge theoretic
argument to higher dimension and showed that the Shafarevich
conjecture holds for any smooth projective variety   with  fundamental
group having a
faithful representation, Zariski dense in a reductive complex
algebraic group - see \cite{Eys}.

Several influential
contributions to these and closely related topics were also made by Lasell
and the fourth author \cite{LasRam}, Mok \cite{M1} and Zuo \cite{Zuo}.

The present article   studies the conjecture in the case when $\pi_1(X,x)$ has
a finite dimensional complex linear representation with infinite
monodromy group.  It combines and develops further  some known techniques in Non abelian Hodge theory.
In particular we
prove the conjecture for projective manifolds $X$ whose fundamental
group admits a faithful representation in $GL_n(\mathbb C)$.

The general strategy has two main steps. First we
use the given faithful linear representation to construct certain
complex variations of mixed Hodge structures ($\mathbb{C}$-VMHS). Then
we utilize the associated period mappings to construct a Shafarevich
morphism. This is quite similar to the way period maps for complex
variations of pure Hodge structures ($\mathbb{C}$-VHS) were used in
\cite{Eys}.
Once the Shafarevich morphism is constructed, holomorphic convexity is
much simpler to obtain here. The crucial point in the construction of the
Shafarevich morphism is a rather subtle rationality lemma which turned
out to rely on Mixed Hodge Theory.

The paper is organized as follows. Section  2 introduces Absolute
Constructible Sets 
  and recalls results from  \cite{Eys}. Section 3
introduces a
$\C$-VMHS constructed in \cite{ES} which serves as a main ingredient of the proof. 
 Section 4 contains the proof of an
important strictness statement. Section 5 contains a rationality lemma and the
reduction to finite number of local systems.   Section 6 contains the
construction of the Shafarevich morphism and the proof of the main
theorem.

 Given present-day technology, it seems difficult to go significantly
 further  in the direction of
proving the  Shafarevich conjecture. Perhaps, the generalization
 to the K\"ahler case or understanding sufficient conditions 
for holomorphic convexity of the universal covering space 
of a singular projective variety might produce interesting developments. 
Several interesting observations have been made in cases of  nonresidually
finite fundamental groups. Bogomolov and the second author suggest
\cite{BK} that the Shafarevich conjecture might fail in the case of
nonresidually finite fundamental groups. From another point of view,  papers by
Bogomolov and de Oliveira \cite{BO1, BO2} suggest that big part of
universal coverings of smooth projective varieties might  still be
holomorphically convex.

\subsection*{Notations}
\addcontentsline{toc}{subsection}{Notations}

 In
  what follows, 
 $X$ will denote a connected projective algebraic complex manifold,
 $x\in X$ a point, $ \Qbar \subset \ell \subset \C$ a field of
 definition for $X$ , and $Z$ a connected projective algebraic
 variety.
 
\subsection*{Statement of the Main Theorem}
\addcontentsline{toc}{subsection}{Statement of the Main Theorem}

\begin{theoi}\label{main} Let $G$ be a reductive algebraic group
defined over $\Q$.  Let $M=M_B(X,G)$ be the character scheme of
$\pi_1(X,x)$ with values in $G$.  
\begin{itemize} 
\item[{\bfseries(a)}] Let 
$\widetilde{H^{\infty}_{M}} \subset
\pi_{1}(X,x)$ be the intersection of the kernels of all
representations $\pi_1(X,x)\to G(A)$, where $A$ is an arbitrary
$\C$-algebra of finite type.  Then, the associated 
Galois covering space of $X$:
\[
\widetilde{X^{\infty}_{M}}= \widetilde{X^{\text{univ}}} /
\widetilde{H^{\infty}_{M}}
\] 
is holomorphically convex.
\item[{\bfseries (b)}]
There exists a natural non-increasing family 
\[
\widetilde{H^{1}_{M}} \subseteq \widetilde{H^{2}_{M}} \subseteq \cdots
\subseteq \widetilde{H^{k}_{M}} \subseteq \cdots \subseteq
\widetilde{H^{\infty}_{M}} \trianglelefteq \pi_{1}(X,x)
\]
of normal subgroups in $\pi_{1}(X,x)$. For a given $k$ the group
$\widetilde{H^{k}_{M}}$ 
corresponds to representations $\pi_1(X,x)\to G(A)$, with $A$ an
Artin local algebra, and such that the Zariski closure of their monodromy
group has $k$-step unipotent radical. For every
$\widetilde{H^{k}_{M}}$  the associated cover
\[
\widetilde{X^k_{M}}=
\widetilde{X^{\text{univ}}} / \widetilde{H^{k}_{M}}
\] 
is
holomorphically convex. 
\end{itemize}
\end{theoi}

{\bfseries Remarks.}
 {\em Let $\pi:\widetilde{X^k_{M}}\to \widetilde{S^k_{M}}(X)$ be the 
Cartan-Remmert reduction of $\widetilde{X^k_{M}}$. The quotient group 
$\pi_1(X,x)^k_{M}:=\pi_1(X,x)/\widetilde{H^k_{M}}$ acts properly discontinuously on $\widetilde{X^k_{M}}$
and $\pi$ is equivariant. We then define the Shafarevich variety as the normal compact complex space
$Sh^k_M(X)= \widetilde{S^k_{M}}(X)/ \pi_1(X,x)^k_{M}$. The resulting Shafarevich morphism
$sh^k_{M}: X \to Sh^k_{M}(X)$ is then independent on $k\in \N^*\cup\infty$.

For every subgroup $H\subset \pi_1(X,x)$ such that $\widetilde{H^{\infty}_{M}}\subset H\subset \widetilde{H^1_{M}}$
the covering space $\widetilde{X^{univ}}/H$ is holomorphically convex as well. }

See section \ref{section-hkm} for the precise definition of $\widetilde{H^k_{M}}$.  If $G=GL_1$, this theorem is a restatement of
\cite{Katnilp}.  Actually, the theorem is likely to hold when we
replace $M_B(X,G)$ by an arbitrary absolutely closed set $M$ defined
over $\Q$ \cite{Simrg1}.

{\bfseries Acknowledgments:}

We thank Fr\'ed\'eric Campana,  J\'anos Koll\'ar and Carlos Simpson for useful
conversations on the Shafarevich conjecture and non-abelian Hodge
Theory.

We have to apologize for the excessive delay between our first
announcement talks on this subject and the availability of a text in
preprint form.

\section[Absolute Constructible Sets]{Absolute Constructible
  Sets}
 
\subsection{Basic facts}
  
Let $G$ be an algebraic reductive group defined over $\Qbar$. The
{\em representation scheme} of $\pi_{1}(X,x)$ is an affine
$\Qbar$-algebraic scheme described by its functor of
points:
$$ 
R(\pi_1(X,x),G)(\mathrm{Spec} (A)) := \mathrm{Hom}(\pi_1(X,x),G(A)) 
$$ 
 for any $\Qbar$ algebra $A$. The {\em character
scheme} of
$\pi_1(X,x)$ with  values in $G$ is the affine scheme $$M_B(X,G)= 
R(\pi_1(X,x),G)//G.$$

Let $\bar k$ be an algebraically closed field of
characteristic zero. Then $M_B(X,G) (\bar k)$ is the set of $G(\bar
k)$-conjugacy classes of {\em reductive} representations of
$\pi_1(X,x)$ with values in $G(\bar k)$, see \cite{LuMa}.

Character schemes of fundamental groups of complex projective
manifolds are rather special.  In \cite{Sim3}, two additional
quasi-projective schemes over $\ell$ are constructed: $M_{DR}(X,G)$
and  $M_{Dol}(X,G)$. The $\C$-points of $M_{DR}(X,G)$ are in
  bijection with the equivalence classes of flat $G$-connections with
  reductive monodromy, and the $\C$-points of $M_{Dol}(X,G)$ are in
  bijection with the isomorphism classes of polystable $G$-Higgs
  $G$-bundles with
  vanishing first and second Chern class. 
Whereas the notion of a
polystable Higgs bundle depends on the choice of a polarization on $X$
 the moduli space $M_{Dol}(X,G)$ does not, i.e. - all moduli
spaces one constructs for the different polarizations are naturally
isomorphic, \cite{Sim3}. This is analogous to the classical
statement that the usual Hodge decomposition on the de Rham cohomology
is purely complex analytic, i.e. independent of a choice of a K\"{a}hler metric\footnote{The harmonic representative
of a cohomology class depends in general on the K\"ahler metric. A helpful remark in the present context is that the harmonic representative 
of a degree $1$ cohomology class actually does not depend on the K\"ahler metric.}. $M_{Dol}(X,G)$ is
acted upon algebraically by the multiplicative group $\C^*$. There is
furthermore a complex analytic biholomorphic map
$$
RH: M_B(X,G)(\C)\to M_{DR}(X,G) (\C)
$$ 
and a real
analytic homeomorphism
$$
KH: M_B(X,G)(\C)\to M_{Dol}(X,G) (\C).
$$ 
$RH$ and $KH$ are also independant of the choice of a K\"ahler metric. 
When $l=\Qbar$, one defines an {\em absolute
constructible subset} of $M_B(X,G)(\C)$
to be a subset $M$ such that: 
\begin{itemize}
 \item $M$ is the set of complex points of a
$\Qbar$-constructible subset of $M_B(X,G)$,
\item $RH(M)$ is the set of
complex points of a $\Qbar$-constructible subset of $M_{DR}(X,G)$,
\item $KH(M)$ is a $\C^*$-invariant set of complex points of a
$\Qbar$-constructible subset of $M_{Dol}(X,G)$.
\end{itemize}

There is a rich theory describing the structure of absolutely
  constructible subsets in $M_B(X,G)$. Here we briefly summarize only those
  properties of absolutely constructible sets that we will need
  later. Full proofs and details can be found in \cite{Simrg1}.

\begin{itemize}
\item The full moduli space $M_B(X,G)$ of representations of $\pi_1(X,x)$
in $G$ defined in \cite{Sim3} is absolutely  constructible and quasi
compact (acqc).
\item The closure (in the classical topology) of an acqc subset is
  also acqc.
\item  Whenever $\rho$ is an isolated point in $M_B(X,G)$, $\{ \rho \}$
is acqc.
\item Absolute constructibility is invariant under standard geometric
    constructions. For instance,
  for any morphism $f : Y \to X$ of smooth connected projective
  varieties, the property of a subset being absolutely constructible
  is preserved when taking images and preimages via $f^{*} :
  M_B(X,G)\to M_B(Y,G)$. Similarly, for any homomorphism $\mu: G\to
  G'$ of reductive groups, taking images and preimages under
  $\mu_*:M_B(X,G)\to M_B(X,G')$ preserves absolute constructibility.
\item Given a dominant morphism $f:Y\to X$ and $i \in \mathbb N$ the
  set $M^i_f(X,GL_n)$ of local systems $V$ on $Y$ such that $R^i f_*V$
  is a local system is ac. Also, taking images and
  inverse images under $R^if_*:M^i_f(X,GL_n)\to M_B(Y,GL_{n'})$
  preserves acqc sets.
\item The complex points of a
 closed acqc set $M$ are stable under the $\mathbb C^*$ action defined
 by \cite{Sim1} in terms of Higgs bundles. By \cite{Sim1} the fixed
 point set $M^{\text{VHS}} := M^{\mathbb{C}^*}$ consists of
 representations underlying
 polarizable complex Variations of Hodge structure ($\C$-VHS, for
 short). Furthermore $M$ is then the smallest closed acqc set in
 $M_{B}(X,G)$ containing
 $M^{\text{VHS}}$.
\end{itemize}

 \subsection{Reductive Shafarevich conjecture}

 After complete results were obtained for surfaces in \cite{KR}, the
 Shafarevich conjecture on
 holomorphic convexity for reductive linear coverings of
 arbitrary projective
 algebraic manifolds over $\mathbb C$ was settled affirmatively in \cite{Eys}.

 \begin{theo}
Let $M\subset M_B(X,G)$ be an absolute constructible set of conjugacy classes of
linear reductive representations of $\pi_1(X,x)$ in some reductive
algebraic group $G$ over $\Qbar$.

Define a normal subgroup  $H_M \subset \pi_1(X,x)$ by:

$$H_M=\bigcap_{\rho \in M(\overline{\mathbb Q})}
\ker (\rho).
$$

The Galois covering space  $\widetilde{X_M}=
\widetilde{X^{\text{univ}}} / H_M$ is holomorphically convex.
\end{theo}

Without a loss of generality we may assume in this theorem that
  $M$ is a closed absolutely constructible set since
  $\widetilde{X_ M}=\widetilde{X_{\bar M}}$.

Let $\Gamma_M$ be the quotient group defined by $$\Gamma_M=\pi_1(X,x) / \bigcap_{\rho \in
  M(\overline{\mathbb Q})} \ker (\rho). $$ $\Gamma_M$ is the Galois group of
  $\widetilde{X_M}$ over $X$. and acts in a proper
discontinuous fashion on the Cartan-Remmert reduction
$\widetilde{S_M}(X)$ of $\widetilde{X_M}$, which is a normal complex
space.
The quotient space $$Sh_M(X)=\widetilde{S_M}(X)/ \Gamma_M$$ is then a normal
projective variety and the quotient morphism $sh_M:X\to Sh_M(X)$ is
called the Shafarevich morphism attached to $M$. This morphism is a fibration, i.e.: 
is surjective with connected fibers. 

Its fibers $Z$ are connected, have the property that $\pi_1(Z)\to
\Gamma_M$ has finite image and are maximal with respect
to these properties. 
\begin{coro} If $\pi_{1}(X,x)$ is almost reductive (i.e. has a Zariski dense representation  with finite kernel in a reductive algebraic group over $\C$)
  then the Shafarevich conjecture holds for $X$.
  \end{coro}

\section[$\C$-VMHS attached to an absolute closed set]{$\C$-VMHS
  attached to an absolute closed set} 

 We will first review some of the results in
\cite{Ha3} and \cite{ES} that enable one to construct various $\C$-VMHS on
$X$ out of $M$.

The results in \cite{Ha3} are important, general and abstract since
they deal with general compactifiable K\"ahler spaces.  The
results in \cite{ES} deal with the less general situation of a compact
K\"ahler manifold but are more explicit and give some useful
properties we have not been able to deduce from
\cite{Ha3}.  More to the point, \cite{ES} will be
sufficient for proving the Shafarevich
conjecture in the case when $\pi_1(X,x)$ has a faithful complex linear
representation. On the other hand, the results in \cite[Sections
  1-12]{Ha3} are needed for the optimal form of our results.

\subsection{$\C$-VMHS, definition, basic properties}

The notion of polarized $\C$-VHS was introduced in \cite{Sim1} as a
straightforward variant of \cite{Gri}.  We will use another equivalent
definition:

\begin{defi}
A $\C$-VHS (polarized complex variation of Hodge structures)
  on $X$ of weight $w\in \Z$ is a 5-tuple $(X,
  \mathbb{V},\mathcal{F}^{\bullet},\overline{ \mathcal {G}}^{\bullet},
  S)$ where:
\begin{enumerate}
\item $\mathbb{V}$ is a local system of finite dimensional
  $\C$-vector spaces, 
\item   $S$ a non degenerate flat sesquilinear pairing on
  $\mathbb{V}$,  
\item $\mathcal{F}^{\bullet} =(\mathcal{F}^p)_{p\in \Z}$ a biregular
  decreasing filtration of $\mathbb{V}\otimes_{\mathbb C}
  \mathcal{O}_X$ 
by locally free coherent analytic sheaves    such that $d' {\mathcal F}^p \subset \mathcal F ^{p-1} \otimes
  \Omega^1_X$, 
\item $ \overline{ \mathcal {G}}^{\bullet} =(\overline{ \mathcal
  {G}}^{q})_{q\in\Z}$ a biregular decreasing filtration of
  $\mathbb{V}\otimes_{\mathbb C} \mathcal{O}_{\bar X}$ 
by locally free coherent antianalytic sheaves such that
 $d''\overline{\mathcal{G}}^p \subset \overline{\mathcal{G}
}^{p-1} \otimes \Omega^1_ {\bar X}$,
\item for every point $x\in X$ the fiber at $x$ $
  (\mathbb{V}_x,\mathcal{F}_x^{\bullet},\overline{ \mathcal {G}}_x^{\bullet})$ is a $\C$-MHS
  polarized by 
  $S_x$.
\end{enumerate}
\end{defi}

The conjugate
$\C$-VHS is the $\C$-VHS obtained on $\overline{\mathbb{V}}$ setting
$\mathcal{F}_{\overline{\mathbb{V}}}^{\bullet} = \overline{ \mathcal
{G}}^{\bullet}$, etc. The local system  $\mathbb{V}\oplus\overline{\mathbb{V}}$
carries a real polarized Variation of Hodge Structures.

Recall that a real reductive algebraic group $E$ is said to be of
Hodge type if there is a morphism $h:U(1)\to Aut(E)$ such that $h(-1)$
is a Cartan involution of $E$, see \cite[p.46]{Sim2}.  By definition,
$h$ is a Hodge structure on $E$.  Connected groups of Hodge type are
precisely those admitting a compact Cartan subgroup.  A Hodge
representation of $E$ is a finite dimensional complex representation
$\alpha: E(\R)\to GL (\mathbb{V}_{\C})$ such that $h$ fixes $\ker
(\alpha)$.In this
case, $\mathbb{V}_{\C}$ inherits a pure polarized Hodge
structures of weight zero. The adjoint representation of a Hodge group
is Hodge. Thus the Lie algebra $\mathfrak{E}$ of $E$ has a natural
real Hodge structure of weight $0$ compatible with the Lie bracket. The
Lie algebra action $\mathfrak{E}_{\C} \otimes \mathbb{V}_{\C} \to
\mathbb{V}_{\C}$ respects the Hodge structures.

The real Zariski closure $E_{\rho}$ of the monodromy group of a
representation $\rho:\pi_1(X,x)\to G(\C)$ underlying a $\C$-VHS is a
group of Hodge type. We have $E_{\rho}\subset R_{\C|\R}
G_{\C}$, where $R_{\C|\R}$ is the Weil restriction of
  scalars functor.
Every Hodge representation $\alpha$ of $E$ gives rise to $$\alpha \circ
\rho: \pi_1(X,x) \to GL( \mathbb{V}_{\alpha})$$ a representation that
underlies a $\C$-VHS \cite[Lemma 5.5]{Sim2}.
 
  The notion of $\C$-VMHS (or graded-polarized
   variation of  $\mathbb{C}$-mixed Hodge structures) used in \cite{ES}
  is a straightforward generalization of that given in
  \cite{SZ, U}:
  
\begin{defi}
A $\C$-VMHS
on $X$ is a 6-tuple $(X,
\mathbb{V},\mathbb{W}_{\bullet},\mathcal{F}^{\bullet},\overline{
\mathcal {G}}^{\bullet}, (S_k)_{k\in\Z})$ where:
\begin{enumerate}
\item $\mathbb{V}$ is a local system of finite dimensional
  $\C$-vector spaces, 
\item $\mathbb{W}_{\bullet}= (\mathbb{W}_k) _{k\in \Z}$ is a
  decreasing filtration of $\mathbb{V}$ by local subsystems,   
\item $\mathcal{F}^{\bullet} =(\mathcal{F}^p)_{p\in \Z}$ a biregular
  decreasing filtration of $\mathbb{V}\otimes_{\mathbb C}
  \mathcal{O}_X$ 
by locally free coherent analytic sheaves  such that $d' {\mathcal F}^p \subset \mathcal F ^{p-1} \otimes
  \Omega^1_X$, 
\item $ \overline{ \mathcal {G}}^{\bullet} =(\overline{ \mathcal
  {G}}^{q})_{q\in\Z}$ a biregular decreasing filtration of
  $\mathbb{V}\otimes_{\mathbb C} \mathcal{O}_{\bar X}$ 
by locally free coherent antianalytic sheaves such that $d''\overline{\mathcal{G}}^p \subset \overline{\mathcal{G}
}^{p-1} \otimes \Omega^1_ {\bar X}$, 
\item  $\forall x\in X$ the stalk $(\V_x,  \mathbb{W}_{\bullet,x},\mathcal{F}^{\bullet}_x,\overline{
\mathcal {G}}^{\bullet}_x)$ is a $\C$-MHS,
\item  $S_k$ is flat sesquilinear non degenerate pairing on
  $Gr^{\mathbb{W}}_k \mathbb{V}$, 
\item $(X, Gr^{\mathbb{W}}_k {\mathbb V},   
  Gr^{\mathbb{W}\otimes_{\C} O_X}_k\mathcal{F}^{\bullet},
   Gr^{\mathbb{W}\otimes_{\C} O_{\bar X}}_k\overline{\mathcal{G}}^{\bullet}, S_k)$ is a $\C$-VHS.  
\end{enumerate}
\end{defi}

We use the following terminology in the sequel:
\begin{defi}
 A homomorphism of groups $\rho: \Gamma \to \Gamma'$ will be called trivial if $\rho(\Gamma)=\{ e\}$. 
A VMHS will be called trivial (or constant) if its monodromy representation is trivial. 
\end{defi}

\subsection{Mixed Hodge theory for the relative completion}

In \cite[Theorem~13.10]{Ha3}, certain $\R$-VMHS are attached to a
$\R$-VHS on compact K\"ahler manifold. In this section we review the results of \cite{Ha3}
  relevant to our discussion and complement them with some explicit
  examples. We will omit the proofs of the statements that
are not essential to our present goals but we will describe in greater
detail the  examples we need.
   
\subsubsection{Hain's theorems}
\label{ssec-relative} 
Let us first review \cite[Sections 1-12]{Ha3}. Another reference where this material (and much more) has been  nicely rewritten in a more general form is \cite[Section 6]{Pri2}.
Let $E^{\rho}$ be a real reductive group of
Hodge type. Let $\rho:\pi_1(X,x)\to E^{\rho}(\R)$ be a Zariski dense
representation underlying a VHS,  and let
 $$1\to \mathcal{U}^{\rho}_x\to \mathcal{G}^{\rho}_x\to E^{\rho}\to 1, 
 \ \ a:\pi_1(X,x)\to \mathcal{G}^{\rho}_x(\mathbb R)$$
 be its relative completion \cite{Ha3}.
 
 $\mathcal G^{\rho}_x$ is a proalgebraic group over $\R$ and $\mathcal
 U^{\rho}_x$ is its prounipotent radical.
 
 Let $k\ge 1$ be an integer, $\mathcal{U}^{\rho}_{x,k}$ be the $k$-th
term of the lower
 central series of  $\mathcal
 U^{\rho}_x$ and $\mathcal{G}^{\rho}_{x,k}$ be
 $\mathcal{G}^{\rho}_x/\mathcal{U}^{\rho}_{k_x}$.
 
 The commutative Hopf algebra $\R[\mathcal{G}^{\rho}_x]$ of the
 regular functions on $\mathcal G^{\rho}_x$ carries a compatible
 $\R$-MHS with nonnegative weights.  The increasing weight filtration
 is described by the formula:
$$ 
W_k \R[\mathcal{G}^{\rho}_x] = \R [\mathcal{G}^{\rho}_{x,k}]
$$ 
where $\R [\mathcal{G}^{\rho}_{x,k}]$ is identified with its image
in $\R [\mathcal{G}^{\rho}_{x}]$. Although these
  MHS are not necessarily finite dimensional, they are always filtered
direct limits of finite dimensional ones.
 
 Let $\mathbb{M}_x=(M_x,W_{\bullet }, F^{\bullet})$ be a finite
 dimensional complex mixed Hodge structure and consider $\alpha:
 \mathcal{G}^{\rho}_x(\mathbb R) \to GL (M_x)$ a representation of
 $\mathcal{G}^{\rho}_x$.  We will say
$\alpha$ is a Mixed Hodge representation iff $\alpha$ is the
 representation arising from the real points of a rational
 representation of $\mathcal G^{\rho}_x$ in $M_x$ and the coaction
$$
\alpha^*: M_x \to M_x \otimes \C[\mathcal{G}^{\rho}_x]
$$ 
respects the
 natural MHS.

 The main result of \cite[Section~13]{Ha3} can now be stated as
   follows.

\begin{prop} \label{prop-relative} Let $\alpha$ be a Mixed Hodge representation.
The representation $\alpha \circ a:\pi_1(X,x) \to GL(M_x)$ underlies a
$\C$-VMHS.  Moreover, any $\C$-VMHS $\mathbb M$ whose graded
constituents $Gr_W ^{k} \mathbb{M}$ are VHS such that their monodromy
representations $\pi_1(X,x) \to GL (Gr_W ^{k} \mathbb{M}_x)$ factor
through $\rho$ is of this type.  
A similar statement holds for $\R$-VMHS.
 \end{prop}

The recent preprint \cite{Ara} gives among other things an alternative approach to this material.

\subsubsection [Example]{Example}
\label{ssec:example}
In \cite{Ha3} Hain describes the steps
$\mathbb{W}_{k}\mathbb{M}$ of the weight filtration in
Proposition~\ref{prop-relative} through iterated integrals. This
however is somewhat technical and goes beyond the scope
of the present paper. 
Instead of discussing the general construction,
we will spell out the definition of the  $\C$-VMHS underlying some
very specific $\C$-Mixed Hodge 
representations of $\mathcal{G}^{\rho}_{x,1}$ which will play a
prominent role in our considerations.
 
Let $(X,\mathbb{V},\mathcal{F}^{\bullet},\overline{ \mathcal {G}}^{\bullet},
S)$ be a $\C$-VHS that will be assumed with no loss of generality of
weight $0$. We will write $\V$
for short, since this will not cause any confusion.

Let $\mathcal{E}^{\bullet}(X,\mathbb V)$ be the de Rham complex
of $\mathbb V$. This de Rham complex inherits a Hodge filtration from
$\mathcal{F}$ and the Hodge filtration on $\mathcal{E}^{\bullet } (X)$
and an anti-Hodge
filtration from $\mathcal{G}$ and the anti-Hodge filtration on
$\mathcal{E}^{\bullet } (X)$. The resulting two filtrations on its
cohomology groups  define on $H^p(X, \mathbb V)$ a $\C$-Hodge
structure of weight $p$.
Furthermore, once we fix a K\"ahler form on $X$, there is a subspace
$\mathcal{H}^p(X, \mathbb V) \subset \mathcal{E}^p(X,\mathbb V)$ consisting of
harmonic forms in a suitable sense such that the composite map $[ \_
]: \mathcal{H}^p(X, \mathbb V) \subset \mathcal{Z}^p(X, \mathbb V)\to
H^p(X, \mathbb V)$ is an isomorphism. This is standard and can be found in e.g. \cite{Zu} for $\mathbb V$ a $\R$-VHS. The general $\C$-VHS case follows in exactly the same way.

\begin{rem}
When $p=1$, the space of harmonic forms is actually independent of
the K\"ahler metric and of the polarization $S$. Furthermore if $Y\to X$ is a morphism $f^*
\mathcal{H}^1(X, \mathbb V)\subset \mathcal{H}^1(Y, \mathbb V)$.
Indeed $\mathcal{H}^1(X, \mathbb V)= \ker(D') \cap\ker(D'')\cap
\mathcal{E}^{1}(X,\mathbb V)$.
\end{rem}

Consider $\alpha \in \mathcal{H}^1(X, \mathbb V)$ such that $[\alpha]$
is of pure Hodge type $(P,Q)$. Then, for all $y\in X$, $\alpha(p) \in
\mathbb{V}_y^{P-1,Q}\otimes \Omega^{1,0} \oplus
\mathbb{V}_y^{P,Q-1}\otimes \Omega^{0,1}$.

Let $(\alpha_i)_{i\in I}$ be a $\C$-basis of $\mathcal{H}^1(X, \mathbb
V)$ such that each $[\alpha_i]$ is of pure Hodge type. Let
$([\alpha_i]^*)_{i\in I}$ be the dual basis of the dual vector space $H^1(X,\mathbb{V})^*$ and
define 
$$
\Omega \in \mathcal{E}^1 (X, \mathbb{V} \otimes
H^1(X,\mathbb{V})^{*})
$$
by the formula:
$$
\Omega= \sum_i \alpha_i \otimes [\alpha_i]^*.
$$ 

Note that $\Omega$ does not depend on the chosen basis. Now we define a new connection on the
vector bundle  underlying the local
system $\mathbb{M}_0=\C\oplus \V^* 
\otimes H^1(X,\V)$ on $X$ a new connection by setting:
$$ d_{\mathbb{M}} = 
 \left(
 \begin{array}{cc}
 d_{\C} & \Omega \\
 0 &  d_{\V^*}\otimes Id _{H^1(X,\V)}
 \end{array}
 \right). 
$$ 
The duality pairings $H^1(X, \V) \otimes H^1(X, \V)^* \to \C$ and $\V
\otimes \V^* \to \C$
are tacitly used in this formula.  Since $d_{\V}\alpha_i=0$, the connection follows
that $d_{\mathbb{M}}$ is a flat connection and this gives rise to a local system $\mathbb{M}$. Furthermore, the connection $d_{\mathbb{M}}$ respects the 2-step
filtration $W^0\M= \C$ $W^1\M=\M$, hence
$\M$  is a filtered
 local system whose graded parts are $Gr_W^0 \M =\C$ and
$Gr_W^1\M= \V^* \otimes H^1(X,\V)$.
 
We now define a Hodge filtration $\mathcal{F}^{\bullet}$ of the
smooth vector bundle underlying $\M$ by the formula valid for every $p\in X$:
$$ 
\mathcal{F}^{k}_p= 
\begin{cases}
\mathcal{F}^{k} _{\V^*_p\otimes H^1(X,\V) }\subset
\V^*_p \otimes H^1(X,\V)\ & \ \mathrm{if} \ \ k>0, \\
\C_p \oplus \mathcal{F}^{k} _{\V^*_p \otimes H^1(X,\V) }\subset \C_p\oplus
\V_p^* \otimes H^1(X,\V)\ & \ \mathrm{if} \ \  k\le 0. 
\end{cases}
$$
Similarly, one defines an anti-Hodge filtration on $\M$ which we denote by
$\overline{ \mathcal {G}}^{\bullet}$ by the formula valid for every $p\in X$:
$$ 
\overline{\mathcal{G}}^{k}_p= 
\begin{cases}
\overline{\mathcal{G}}^{k} _{\V^*_p\otimes H^1(X,\V) }\subset
\V^*_p \otimes H^1(X,\V)\ & \ \mathrm{if} \ \ k>0, \\
\C_p \oplus \overline{\mathcal{G}}^{k} _{\V^*_p \otimes H^1(X,\V) }\subset \C_p\oplus
\V_p^* \otimes H^1(X,\V)\ & \ \mathrm{if} \ \  k\le 0. 
\end{cases}
$$

It defines on each stalk 
$\M_p$ a $\C$-MHS such that $Gr^0_W \M_p$ is the trivial Hodge structure on $\C$ and $Gr^1_W \M_p$ is the given Hodge Structure
on $\V_p^* \otimes H^1(X,\V)$.

\begin{lem}
$\mathcal{F}^{k}$ is a holomorphic subbundle of the holomorphic
vector bundle $\mathcal{M}$ underlying $\M$
 and
 satisfies Griffiths transversality. 
\end{lem}
\begin{pf}
First observe that
the $\bar\partial$ operator of $\mathcal{M}$ is given  by $d_{\M}^{0,1}$.
Consider the original flat connection

$$d = 
 \left(
 \begin{array}{cc}
 d_{\C} & 0\\
 0 &  d_{\V^*}\otimes Id _{H^1(X,\V)}
 \end{array}
 \right).
$$ 
Obviously 
$$
d_{\M}^{0,1} = d^{0,1}+ \left(
 \begin{array}{cc}
 0 & \Omega^{0,1} \\
 0 & 0
 \end{array}
 \right).
$$ 
Since $d^{0,1}$ preserves $\mathcal{F}^{k}$, it follows that $\mathcal{F}_k$ is an
holomorphic
subbundle of $\mathcal{M}$ iff $$\left(
\begin{array}{cc}
 0 & \Omega^{0,1} \\ 0 & 0
 \end{array}
 \right) \mathcal{F}^{k} \subset \Omega^{0,1} \otimes
 \mathcal{F}^{k},$$ where $\Omega^{0,1} \in
 \mathcal{E}^{0,1}(X,\V)\otimes H^1(X,\mathbb{V})^*$ 
 is the $(0,1)$-component of $\Omega$. 
This  condition is equivalent to 
 $\Omega^{0,1}\cdot \mathcal{F}^{k} _{\V^* \otimes H^1(X,\V) }=0$ if $k>0$.
 
We are thus reduced to checking that for every $\alpha \in
 \mathcal{H}^1(X, \mathbb V)$ such that $[\alpha]$ is of pure Hodge
 type $(P,Q)$ and $[\beta]^{\vee} \in (H^1(X,\V) ^*)^{-P,-Q}$ $$\alpha^{0,1}
 \otimes [\beta]^{\vee}\cdot \mathcal{F}^{k} _{\V^* \otimes H^1(X,\V) }=0, \ 
 \text{if} \  k>0.$$
 
 It is enough to check that $[\alpha] \otimes [\beta]^{\vee} \cdot
H^{k,1-k} _{\V^* \otimes H^1(X,\V) }=0$, or further decomposing in
Hodge type that 
$$\alpha^{0,1} \otimes [\beta]^{\vee}\cdot   h^{-P'+k,
-Q'-k+1} \otimes [\beta] ^{ P',Q'}=O$$
 where $ h^{-P'+k, -Q'-k+1}\in
(\V^*)^{-P'+k, -Q'-k+1}$ and $ [\beta] ^{ P',Q'}\in H^1(X,\V)^{
P',Q'}$. The only non trivial case is when $P'=P, Q'=Q$ and this
reduces to showing that $ \mathbb{V}^{P,Q-1}\otimes \Omega^{0,1}
. (\V^*)^{-P+k, -Q-k+1}=0$ which is the case since $k>0$.
 
Griffiths transversality is the statement that $
d^{1,0}_{\M}\mathcal{F}^{k}\subset \mathcal{F}^{k-1} \otimes
\Omega^{1,0}$ and follows from the same argument.
\end{pf}

Antiholomorphicity and Griffiths anti-transversality for $\overline{\mathcal{G}^{\bullet}}$ 
can be proved by the same method. Hence we have defined on $\M$ a graded polarizable VMHS with weights
$0,1$, the polarizations being the natural ones. In \cite{HZ}, the
case of $\V=\C_X$ is treated. In that case, the VMHS is actually
defined over $\Z$.

\begin{defi}
$\M=\M(\V):=(X,
\mathbb{M},\mathbb{W}_{\bullet},\mathcal{F}_M^{\bullet},\overline{
\mathcal {G}}_M^{\bullet}, (S_k)_{k=0,1})$ is the $1$-step $\C$-VMHS
attached to $\V$.
 \end{defi}

\subsection{Mixed Hodge theory for the deformation functor}\label{mhdef}

In this paragraph, we review the construction of \cite{ES}. The \lq new\rq \  aspects of this construction 
actually grew out of the previous example. 
The older aspects, on the other hand, were part of Goldman-Millson's theory of deformations for representations of K\"ahler groups \cite{GM}.

In this paragraph, we fix $N\in \N$ and  assume that $G=GL_N$ and $M=M_B(X, GL_N)$.
Let $\rho:\pi_1(X,x)\to GL_N(\C)$ be the monodromy represention of a $\C$-VHS.
Let $\hat\Or$ be the complete local ring of $[\rho] \in\rep(\C)$. Let 
$$\text{\sf obs}_{2}=[-;-]: S^2 H^1(X, \mathrm{End}(\V_{\rho}))\to H^1(X, \mathrm{End}(\V_{\rho}))
$$
be the Goldman-Millson obstruction to deforming $\rho$. Define $I_2, (I_n)_{n\ge 2}, (\Pi_n)_{n\ge 0}$, as follows: 
\begin{eqnarray*}
\Pi_0 &=& \C \\
\Pi_1 &=& H^1(X, \mathrm{End}(\V_{\rho}))^* \\
I_2&=&\mathrm{Im}(^t\text{\sf obs}_{2}) \subset S^2 H^1(X, \mathrm{End}(\V_{\rho}))^* \\
I_n&=& I_2 S^{n-2} H^1(X, \mathrm{End}(\V_{\rho}))^*\\
\Pi_n&=& S^n H^1(X, \mathrm{End}(\V_{\rho}))^*/I_n
\end{eqnarray*}

Then the complete local $\C$-algebra $$(\hat{\mathcal{O}}_T,\m):=
(\sum_{n\ge 0} \Pi_n, \sum_{n\ge 1} \Pi_n)$$
is the function algebra of a formal scheme $T$ which is the germ at $0$ of the quadratic cone 
$$ \text{\sf obs}_{2}^{-1}(0)\subset H^1(X, \mathrm{End}(\V_{\rho})).$$ 

We endow $\hat{\mathcal{O}}_T$ with a split mixed Hodge structure with non positive weights, whose weight filtration is given by the formula $W_k  \hat{\mathcal{O}}_T = \m^{-k}$ for $k\le 0$,  arising from the identifications:
$$\hat{\mathcal{O}}_T = \sum_{n\ge 0} \m^n/ \m^{n+1} =\sum_{n\in \N} \Pi_n,
$$
$\Pi_n$ being endowed with its natural  $\C$-Hodge structure of weight $-n$. This mixed Hodge structure is infinite dimensional, but can be described as the limit 
of the resulting finite dimensional MHS on $\hat{\mathcal{O}}_T/ \mathfrak{m}^n$.

In \cite{GM}, an isomorphism between $Spf(\hat\Or)$ and $T\times A$ is constructed, 
 where $A$ is the germ at zero of a finite dimensionnal vector space. In \cite{ES}, 
this constuction is revisited. A slight reinterpretation of Goldman-Millson theory 
is that one can realize the formal local scheme $T$ as a hull of the deformation functor for $\rho$.
Actually, there are {\em three} preferred such realizations $\mathcal{GM}^c, \mathcal{GM}',\mathcal{GM}''$
which are given by
three canonical representations:
\[
\begin{split}
\rho^{GM}_{T} : & \pi_1(X,x)\to
 GL_N(\widehat{\mathcal{O}}_T) \\
\rho^{GM'}_{T} : & \pi_1(X,x)\to
 GL_N(\widehat{\mathcal{O}}_T)\\
\rho^{GM''}_{T} : & \pi_1(X,x)\to
 GL_N(\widehat{\mathcal{O}}_T)
\end{split}
\]

These three representations are conjugate up to an isomorphism of $T$. 

We can now summarize the results developped by \cite{ES} in the form we shall need:

\begin{defi} Let $\eta_1, \ldots , \eta_b \in E^{\hdot}(X,\mathrm{End}(\V_{\rho}))$ form a basis of the subspace $\mathcal{H}^1(X,\mathrm{End}(\V_{\rho}))$ of harmonic 
twisted one forms, each $\eta_i$ being of pure Hodge type $(P_i,Q_i)$ for the Deligne-Zucker
$\C$-Mixed Hodge Complex $E^{\hdot}(X,\mathrm{End}(\V_{\rho}))$. Then  $\{ \eta_i \}$ is a basis of
$H^{1}(X,End(\V_{\rho}))$ whose dual basis we denote by $(\{\eta_1\}^*, \ldots, \{\eta_b\}^*)$. 

The $\mathrm{End}(\V_{\rho})\otimes\Pi_1 $-valued one-form   $\alpha^v_1$ is defined by 
the formula:

$$\alpha^v_1=\sum_{i=1}^b \eta_i \otimes \{\eta_i\}^*.
$$
 \end{defi}

\begin{prop} \label{constr}
 For $k\ge 2$, we can construct a unique $D''$-exact form  $\alpha^v_k \in E^1(X,\mathrm{End}(\V_{\rho}))\otimes \Pi_k$  such that the following relation holds:
$$
D'\alpha^v_k +\alpha^v_{k-1}\alpha^v_1 + \alpha^v_{k-2}\alpha^v_2 +\ldots+ \alpha^v_1 \alpha^v_{k-1}=0.
$$
\end{prop}
\begin{prop} Let $A_v=\sum \alpha_k^v$ acting on the vector bundle underlying
the filtered local system $(\V_{\rho} \otimes \hat\Ot, W_k(\V_{\rho}\otimes \hat\Ot)= \V_{\rho}\otimes \m^{k-wght(\V_{\rho})}) $, whose connection will be denoted by $D$.
 
Then, $D+A^v$  respects this weight filtration, satisfies Griffiths' transversality for the Hodge filtration  $\mathcal{F}^{\hdot}$ defined by

$$\mathcal{F}^p=\bigoplus_{k=-n}^0 \mathcal{F}^p(\V_{\rho} \otimes \Pi_{-k})$$
and we can construct an anti-Hodge filtration so that the resulting structure is a graded polarizable $\C$-VMHS whose monodromy representation is $\rho^{GM''}_T$. 

\end{prop}

A detailed proof of this
  proposition is given in \cite{ES}. The essential part is the
  construction of the anti-Hodge filtration which is similar in spirit
  but somewhat subtler than construction given in Example~\ref{ssec:example}.
  
\begin{defi}\label{defi:Dk}
The $\C$-VMHS obtained by reduction $\mod \ \m^n$ corresponds to

\[
\rho_{T,n}:=(\rho^{cGM^v}_{T} \!\! \mod \m^n) : \pi_1(X,x)\to
GL_N(\widehat{\mathcal{O}}_{T}/\mathfrak{m}^{n}).
\]
It  will be called the $n$-th deformation of $\V_{\rho}$ and will be denoted by $\D_n(\V_{\rho})$.
\end{defi}

By construction, $D+A^v$ is an $\hat{\mathcal{O}}_T$-linear connection.
As a consequence of the methods in \cite[pp 18-23]{ES}, we also have:

\begin{prop}\label{mhslice}
There is a MHS on $\hat{\mathcal{O}}_T$ whose weight filtration is given by the
powers of the maximal ideal and such that the natural map $\hat{\mathcal{O}}_T \to \mathrm{End}_{\C}(\D_n(\V_{\rho}))$
respects the natural MHS. 
\end{prop}

This MHS is {\em not} the split MHS constructed above. This split MHS is just the weight graded
counterpart of the true object.  These MHS and VMHS are not 
uniquely defined when the deformation functor of $\rho$ is not prorepresentable.  This
phenomenon does not occur when the representation is irreducible. 

\begin{rem}
 The restriction $G=GL_n$ in the above 
 considerations was introduced only for convenience. It is not essential. In \cite{ES},
 similar statements are proven for arbitrary reductive groups $G$.
\end{rem}

\section{Subgroups of $\pi_1(X,x)$ attached to $M$}

  Let  $G$ be a reductive algebraic group defined over
  $\overline{\mathbb{Q}}$.
Suppose as before $M\subset M_B(X,G)$ is an absolute closed subset. 
\subsection{Definitions}
\label{section-hkm}

\begin{defi}Let $M^{VHS}$ be the subset of $M(\C)$ consisting of  the conjugacy classes of $\C$-VHS that is
$M^{VHS}:= KH^{-1}( M_{Dol}(X, G)^{\C^*}(\C))$. 
\end{defi}

We choose  a set $M^*$ of reductive representations $\rho: \pi_1(X,x) \to G(\C)$ mapping onto $M^{VHS}$ under the natural map $R(\pi_1(X,x), G)\to M_B(X,G)$. 
To be more precise, we define  $M^{*}$ to be the union of the closed
  $G$-orbits on $R(\pi_1(X,x),G)$ -or equivalently the set of reductive complex representations-  whose equivalence class lie in
  $M^{VHS}$. Similarly, we define $M'$  to be the union of the closed
  $G$-orbits on $R(\pi_1(X,x),G)$  whose equivalence class lie in
  $M$.
  To each $\rho\in M^*$ we attach $E^{\rho}$ the real Zariski closure of its monodromy
 group and the other constructions reviewed in paragraph \ref{ssec-relative}.
 
 \begin{defi}  The tannakian categories $\mathcal{T}^{VHS}_M$ and  $\mathcal{T}_M$ are defined as follows:
\begin{description}
\item[$\mathcal{T}^{VHS}_M$]  is the full
  Tannakian subcategory of the category of local systems on $X$
  generated by the elements of $M^{VHS}$.
\item[$\mathcal{T}_M$] is the  full Tannakian subcategory of the
  category of local systems on $X$ 
  generated by the elements of $M$.
\end{description}
 \end{defi}

 Every object in $\mathcal{T}^{VHS}_M$ is isomorphic to an object
which is a  subquotient of $\alpha_1(\rho_1)\otimes\ldots\otimes
 \alpha_s(\rho_s)$, where $\rho_1, \ldots,\rho_s$
 are elements of $M^*$ and
  $\alpha_i$ is a complex
 linear  finite dimensionnal representation of $E^{\rho_i}(\R)$.
Let $M^{**}$ be the set of all such subquotients. 
The objects of  $\mathcal{T}^{VHS}_M$  underly
 polarizable $\C$-VHS.

Let $T^{VMHS}_M$ be the thick Tannakian subcategory of ($\C$-VMHS) whose graded
constituents are
objects of  $\mathcal{T}^{VHS}_M$.
  The full subcategory
of $\mathcal{T}^{VMHS}_M$  with a weight filtration
of length at most  $k+1$ will be denoted by $T^{VMHS}_M(k)$ .

\begin{exe}
 For every $\rho \in M^{**}$, $\alpha$ as above and $\sigma=\alpha\circ \rho$, $\D_k(\V_{\sigma})$ is an object 
of $T^{VMHS}_M(k)$.
\end{exe}

 \begin{defi}
 
 Given $X$, $G$, and $M \subset M_{B}(X,G)$
  as above, and $k \in \mathbb{N}$ we define the following natural
  quotients of $\pi_{1}(X,x)$:
\begin{description}
\item[$\Gamma_M^{\infty}$] is the quotient of $\pi_1(X,x)$ by the
  intersection $H_M^{\infty}$ of the kernels of the objects of
  $\mathcal{T}^{VMHS}_M$ and of the objects of $M$.  
\item[$\widetilde{\Gamma_M^{\infty}}$] is the quotient of $\pi_1(X,x)$
  by the intersection $\widetilde{H_M^{\infty}}$ of the kernels of the
  monodromy representation of $\D_n(\V_{\sigma})$, $\sigma\in M^{**}$,
  $n\in\N$, and of the objects of $M$.
\item[$\Gamma_M^{k}$] is the quotient of $\pi_1(X,x)$ by
the intersection $H_M^{k}$ of the kernels of the objects of
$\mathcal{T}^{VMHS}_M(k)$ and of the objects of $M$.
\item[$\widetilde{\Gamma_M^{k}}$] is the quotient of
$\pi_1(X,x)$ by the intersection $\widetilde{H_M^{k}}$ of the kernels of
the monodromy representation of $\D_k(\V_{\sigma})$, $\sigma\in
M^{**}$, and of the objects of $M$.
\end{description}
 \end{defi}

It is likely that the canonical quotient morphism $\Gamma^k_M \to \widetilde{\Gamma^k_M}$ is an isomorphism
but we do not have a proof of this fact yet. We will thus have to work with the above slightly clumsy notation.

 Note that we have the inclusions:
 
 \[
\begin{split}
\Gamma_M^{\infty} & = \bigcap_{k\in\N}\displaylimits
\Gamma_M^k\subset\Gamma_M^{k+1}\subset \Gamma_M^k
\subset\Gamma^0_M=\Gamma_M, \\
\widetilde{\Gamma_M^{\infty}} & = \bigcap_{k\in\N}\displaylimits
\widetilde{\Gamma_M^k}\subset\widetilde{\Gamma_M^{k+1}}\subset
\widetilde{\Gamma_M^1 }
\subset\widetilde{\Gamma^0_M}=\widetilde{\Gamma_M}.
\end{split}
\]
 It should be noted that since $H^k_M$ (respectively $\widetilde{H^k_M}$) is normal the various base point changing isomorphisms $\pi_X(X,x')\to \pi_1(X,x)$
 respect $H^k_M$ (respectively $\widetilde{H^k_M}$). Hence, dropping the base point dependance in the notation $H^k_M$ (respectively $\widetilde{H^k_M}$) is harmless.

 For future reference, we state the following lemma whose proof is tautological. 
 
\begin{lem} $H^k_M$ is the intersection of $\Gamma_M$ and the kernels of $a^{\rho}_k: \pi_1(X,x) \to \mathcal{G}^{\rho}_{x,k}(\R)$. 
\end{lem}

\subsection{Strictness }

Let $z\in Z$ be a base point in the connected projective variety $Z$. 

\begin{prop}\label{strict}
For every $f:(Z,z)\to (X,x)$
  such that $\pi_1(Z,z)\to \Gamma_M$ is trivial, 
  the following are equivalent:
  \begin{enumerate}
  \item For every $\mathbb V$ in $\mathcal{T}^{VHS}_M$, 
  $H^1(X,\mathbb {V})\to H^1(Z,\mathbb {V})$ is trivial,
  \item $\pi_1(Z,z)\to \Gamma_M^1$ is trivial,
\item $\pi_1(Z,z)\to \widetilde{\Gamma_M^1}$ is trivial,
 \item For every $\mathbb V$ in $\mathcal{T}^{VHS}_M$, 
  for every $\widehat Z_i \to Z$ a resolution of singularities of an irreducible component,
  the VMHS
  $M(\V)_{\widehat Z_i}$ is trivial.
\item For every $\sigma\in M^{**}$ and $k\in \N$ 
, 
  for every $\widehat Z_i\to Z$ a resolution of singularities of an irreducible component,
  the VMHS
  $\D_k(\V_{\sigma})_{\widehat Z_i}$ is trivial.
\item $\pi_1(Z,z)\to \widetilde{\Gamma_M^{\infty}}$ is trivial.
\item $\pi_1(Z,z)\to \Gamma_M^{\infty}$ is trivial.
 \end{enumerate}
\end{prop}

\begin{pf}

$(1\Longleftrightarrow 2)$. 
Fix $\rho\in M^{**}$. Denote by $E^{\rho}$ the real Zariski  closure of $\rho(\pi_1(X,x))$. 
By hypothesis $\rho(\pi_1(Z,z))=\{ e\}$ and thanks to \cite{Ha3}, section 11,  \footnote{Stricto sensu, in order to apply  \cite{Ha3} sect. 11,  we need that $Z$ be a smooth connected manifold.  We can replace $Z$ by a neighborhood $U$ of it in some embedding in $\mathbf{P}^N(\C)$ such that $Z\to U$
is an homotopy equivalence and apply \cite{Ha3} sect. 11 to $U$.} we have a diagram:
$$ \begin{array}{ccccc}
\pi_1(Z,z) & \buildrel{a_Z}\over{\mapsto} & \hat{\pi}^{DR}_1(Z,z)& &  \\
\downarrow & & \downarrow & &\\
\pi_1(X,x) & \buildrel{a_X}\over{\mapsto} & \mathcal{U}^{\rho}_x& \subset& \mathcal{G}^{\rho}_x
   \end{array}
$$
where $\pi_1(Z,z)  \buildrel{a_Z}\over{\mapsto}  \hat{\pi}^{DR}_1(Z,z)= \mathcal{U}^{e}(Z,z)= \mathcal{G}^{e}(Z,z)$ is the Malcev completion of $\pi_1(Z,z)$, i.e.: its relative completion with respect to the trivial
representation. 

Let $\{V_{\alpha}\}_{\alpha}$ be a set of representatives of all
isomorphism classes of complex irreducible left $E^{\rho}$-modules.

The prounipotent group morphism $f_*:\hat{\pi}^{DR}_1(Z,z)\to \mathcal{U}^{\rho}_x $  gives rise to a morphism
of proalgebraic complex vector groups (=limits of finite dimensional
complex vector spaces viewed as algebraic groups):
$$ H_1(\hat{\pi}^{DR}_1(Z,z))( \C) \to H_1 (\mathcal{U}^{\rho}_x)(\C)
$$
where $H_1(\mathcal{U})= \mathcal{U} / \mathcal{U}'$ is the abelianization. One has identifications ( see \cite{Ha3} p 73):

\[
\begin{split}
H_1(\hat{\pi}^{DR}_1(Z,z))( \C)= H_1(Z,\C),\\
H_1 (\mathcal{U}^{\rho}_x)( \C) =\prod_{\alpha} H_1(X,\V_{\alpha}) \otimes V^*_{\alpha}
\end{split}
\]

where $\V_{\alpha}$ is the local system attached to $\rho$ and $V_{\alpha}$. The map is the transpose
of the map $$\bigoplus_{\alpha}  H^1(X,\V^*_{\alpha}) \otimes V_{\alpha} \to H^1(X,\C)$$ given on each factor by the composition:
$$
H^1(X,\V^*_{\alpha}) \otimes V_{\alpha}
 \buildrel{i^*_Z \otimes {\mathrm{id}}_{V_{\alpha}}}\over{\to} H^1(Z,\V^*_{\alpha}) \otimes V_{\alpha}
= H^1(Z,\C)\otimes V^*_{\alpha} \otimes V_{\alpha} \buildrel{{\mathrm{id}}\otimes {\mathrm{tr}}}\over{\to} H^1(Z,\C).$$

For the middle equality in this formula, we used that $\V_{\alpha|Z}$
is the trivial local system, which follows from the assumption that 
$\pi_{1}(Z,z) \to \Gamma_{M}$ is trivial.

Hence condition 1 is equivalent to $ H_1(\hat{\pi}^{DR}_1(Z,z)( \C)) \to H_1 (\mathcal{U}^{\rho}_x)(\C) $ being zero which in turn is equivalent to 
condition 2. 

$(2\Longrightarrow 3)$
Condition 3 is obviously implied by condition 2.

$(3 \Longrightarrow 1)$ If $3$ holds $\D_1(\V_{\sigma})|_Z$ is a trivial local system . But, by construction this local system is a deformation of a 
trivial local system by a one-step nilpotent matrix of closed one forms written in the following block form:
$$\left(\begin{array}{cc}
   0 & A \\
0 &0
  \end{array}\right)
$$ 
hence, in the same basis,  its monodromy on any $\gamma\in\pi_1(Z,z)$ is given by:
$$\left(\begin{array}{cc}
   1& \int_{\gamma} A \\
0 &1
  \end{array}\right).
$$ 
Hence the triviality of $\D_1(\V_{\sigma})|_Z$ implies that $\int_{\gamma} A=0$, or that the cohomology class of $A$ is zero. 
But by construction, the cohomology class of $A$ is zero iff condition $1$ holds.

$(1 \Longrightarrow 4)$ The cohomology class of the form $\alpha_1$ vanishes after restriction to $Z$ and so vanishes after pullback to $\widehat Z_i$. We denote by $f_i$ the composition of  $f$ with the map $\widehat Z_i\to Z$.  But $\alpha_1\in \ker(D')\cap \ker(D'')$.  Hence $f_i^*\alpha_1 \in  \ker(D')_{i,1}\cap \ker(D'')_{i,1}$ where $D'_{i,1},D''_{i,1}$
are the usual $D',D''$  acting on  $E^{1}(\widehat Z_i,\mathrm{End}(f^*_i\V_{\rho}))$.
As mentioned before, Hodge theory implies that: $$ \ker(D')_{i,1}\cap \ker(D'')_{i,1}= \mathcal{H}^{1}(\widehat Z_i,\mathrm{End}(f^*_i\V_{\rho})).$$ 
Hence $f^*_i \alpha_1$ is the harmonic representative of its class. From this it follows that $f_i^*\alpha_i=0$. 
This implies that $f^*_i\M$ is the trivial deformation of $f^*_i\M_0$ and condition 4 follows.

$(4 \Longrightarrow 1)$ The method we used to prove $(3 \Longrightarrow 1)$ works to yield that $H^1(X,\mathbb {V})\to H^1(\hat Z_i,\mathbb {V})$ is zero. But this implies by the argument we used to show $(1 \Longrightarrow 4)$ that $f^*_i\alpha_1=0$. 
This in turn implies that $i_A^* \alpha_1=0$ if $f(Z)=\coprod A$ is a smooth stratification. Hence the holonomy of 
$M(\V)_{Z}$ is trivial. Applying once again the method for $(3 \Longrightarrow 1)$ completes the argument. 

$(1\Longrightarrow 5)$ Continuing the same line of reasoning as in 
proving
$(1 \Longrightarrow 4)$  and using the fact that the $(\alpha^v_k)$ constructed in  proposition \ref{constr}
are uniquely determined, it follows that $(f_i^*\alpha^v_k)$ is the family of twisted forms 
one gets from applying the construction of proposition \ref{constr} starting with $f^*_i \alpha_1=0$. Hence 
$f^*_i \alpha_k^v=0$ and $f^*_iA^v=0$. Condition 5 then follows. 

$(1\Longrightarrow 6)$ Continuing this line of reasoning, the argument made in 
$(4 \Longrightarrow 1)$ implies that the restriction of $\D_k(\V_{\sigma})$ to $Z$ has trivial monodromy, which is equivalent to condition 6. 

$(6\Longrightarrow 2)$ is trivial.

 $(6\Longrightarrow 5)$ comes from the fact that condition 6 implies
that the restriction of $\D_k(\V_{\sigma})$ to $Z$ has trivial monodromy, condition 5 follows a fortiori.

$(7\Longrightarrow 3)$ is trivial. 

$(1 \Longrightarrow 7)$
 The proof is an easy adaptation of the argument of \cite{Katnilp}, section 2. We nevertheless feel it is necessary to give some
 details. 

The Lie algebras $L(Z,z)= {\mathrm{Lie}}(\hat{\pi}^{DR}_1(Z,z))$ and $\mathfrak{U}^{\rho}_x={\mathrm{Lie}}(\mathcal{U}^{\rho}_x)  $
are nilpotent and so come equiped   with a decreasing filtration given by their lower central series. The map $i_{Z}$ gives rise to a Lie algebra
morphism $(i_Z)_*: L(Z,z) \to \mathfrak{U}^{\rho}_x$. It is enough to show that $(i_Z)_*=0$

By relabelling we can convert the lower central series  into an increasing filtration $B^{\bullet}L(Z,z) $ and $B^{\bullet} \mathfrak{U}^{\rho}_x$ with
indices $\le -1$. For both Lie algebras $Gr_B^{-1}(\_)= H_1(\_)$ and $Gr_B^{-1}(L(Z,z))$ generates 
the graded Lie algebra $Gr_B^{\bullet}(L(Z,z))$. Hence condition 3 implies that $Gr_B^{\bullet} (i_Z)_*:Gr_B^{\bullet} L(Z,z) \to Gr_B^{\bullet} \mathfrak{U}^{\rho}_x$
is zero. 

First consider the case where $Z$ is smooth. Then, by \cite{Ha2} \cite{Ha3}, both $L(Z,z)$ and
 $\mathfrak{U}^{\rho}_x$ carry a functorial  Mixed Hodge structure
whose weight filtration is $B^{\bullet}$. Hence, since the map $(i_Z)_*$ 
respects the Mixed Hodge structures, it is strict for the weight filtration and 
$Gr_B^{\bullet} (i_Z)_*=0 \Rightarrow (i_Z)_*=0$.

Next we consider the case where $H_1(Z)$ is pure of weight one. We recall, see \cite{Ha2}, that 
$\R[\widehat{\pi}^{DR}_1(Z,z)]=H^0 (\overline{\text{\sf B}} ( \R,
E^{\bullet}(Z),\R))$ where $\overline{\text{\sf B}}$
 is the reduced bar construction and $E^{\bullet}(Z)$ is a multiplicative mixed Hodge complex computing $H^{\bullet}(Z)$ endowed with a base point at $z$. 
$\overline{\text{\sf B}} ( \R,
E^{\bullet}(Z),\R)$ carries an increasing filtration $\mathfrak{B}_{\bullet}$, the bar filtration.
It follows from \cite{Ha2} that $ \overline{\text{\sf B}} ( \R,
E^{\bullet}(Z),\R)$ endowed with the bar filtration is a filtered mixed Hodge complex so that the bar filtration
on $\R[\widehat{\pi}^{DR}_1(Z,z)]$ is a filtration by MHS. The
 Eilenberg-Moore spectral sequence which is 
the spectral sequence associated to the bar filtration is  a spectral sequence in the category of MHS 
and, since $H^1(Z)$ is pure of weight one,  $E_1^{s,-s}= H^1(Z)^{\otimes s}$ is pure of weight $s$. 
Hence $Gr_{\mathfrak{B}}^k\R[\widehat{\pi}^{DR}_1(Z,z)]$ is pure of weight $k$.  Since the bar filtration is a refinement of the weight filtration,
it follows that the bar filtration and the weight filtration coincide. Combining this with the preceding argument, one easily finishes the
  proof of the case when $H_1(Z)$ is pure of weight one.

  Finally note that by passing to a hyperplane section we may assume
  that $Z$ is a curve which without a loss of generality can be taken
  to be seminormal and the argument of \cite{Katnilp} p. 340-341
applies verbatim. One concludes using lemma 2.4 p. 342 in \cite{Katnilp}.

\end{pf}

\begin{rem}
 If we skip items 2 and 7 in the previous proposition  we obtain a strictness statement which can be proved without relying on 
 \cite{Ha3}.
\end{rem}

\begin{rem}
As far as the equivalence of condition 7 with the rest is concerned, we believe that one can adapt the explicit argument made for $(1\Longrightarrow 6)$ using the more sophisticated iterated
integrals of \cite{Ha3}.  

Except perhaps for condition 7, that depends on $X$ being projective, the proposition is valid in the compact K\"ahler case. 
 \end{rem}
 
 \begin{rem}
A generalization to the K\"ahler case 
of the main result in \cite{Katnilp} with an alternative proof has been given
in the unpublished thesis \cite{Ler} (see also \cite{Cla}) as a byproduct of her exegesis of \cite{Ha2} and \cite{Haalb}. 
The core of her argument could be
reformulated in such a way that it becomes equivalent to the special case of the present one where $G=\{ e \}$ is the trivial group.
\end{rem}

\subsection{Reduction to using VSHM}

\begin{prop} Let $n$ be a non negative integer.
Let $H_n$ be
the
intersection
of the kernels
of all linear
representations
$\pi_1(X)\to GL_n(A)$, $A$ being an arbitrary $\C$-algebra. Let $M=M(X,GL_n)$. 
Then,  $H_n=\widetilde{H^{\infty}_M}$.
\end{prop}

\begin{pf}
 The inclusion $H_n\subset \widetilde{H_M^{\infty}}$ is obvious. Now let $\gamma\in \widetilde{H^{\infty}_M}$. Then $\gamma$ defines a
matrix valued regular function $F$ on $R(\pi_1(X,x), GL_n)$ (i.e.: $F\in \mathrm{Mat}_{n\times n}(\C [R(\pi_1(X,x), GL_n)])$) which reduces to 
the constant function with value $I_n$ on $T_{\rho}\subset R(\pi_1(X,x), GL_n)$ for every element $ \rho\in M^{VHS }$. Goldman-Millson theory implies that the tautological representation $\pi_1(X,x)\to GL_n(\widehat\Or)$ is conjugate to the pull back by $cGM$ of $\rho^{cGM}_T$. Hence $F$ induces the trivial matrix valued function when reduced to $Spf(\widehat\Or)$. Hence $F$ induces the constant matrix valued function with value $I_n$ on some complex analytic neighborhood of $M^{VHS}$. 

Let $\widetilde \rho$ be a semisimple complex representation mapping to $M-M^{VHS}$. Then, by \cite{Sim1}, $\widetilde \rho$ correspond to a polystable Higgs bundle $(\mathcal{E},\theta)$. 
For $t\in\C^*$, let $\tilde \rho (t)$ corresponds to $(\mathcal{E},t.\theta)$. By applying  the Goldman-Millson construction
to each  $\tilde{\rho(t)}$, we get a real analytic family of flat connections $(D_t)_{t\in\C^*}$ on the smooth vector 
bundle underlying $\mathcal{E} \otimes \mathcal{O}_{\tilde\rho(t)}$ (see for instance \cite[pp. 21]{Pri4}) such that the image $F_t$ of the matrix function $F$  in the complete local ring at $\tilde{\rho}(t)$, satisfies $F_t=\mathrm{hol}(D_t)\in \mathrm{Mat}_{n\times n}( \mathcal{O}_{\tilde\rho(t)})$. Since $F_t=I_n$ for small $t$ then $F_1=I_n$. Hence $F$ maps to $I_n$ in $\mathrm{Mat}_{n\times n}(\C [(\Gamma, GL_n)]_{\tilde\rho})$. Hence $F=I_n$ in a complex analytic neighborhood of the
set of semisimple representations. 

Given a non semi simple representation $\rho^{arb}$ we may find a sequence
$(\rho_m)_{m\in \N}$ of  conjugate representations converging to a semi simple one we see that $\rho_m(\gamma)=Id_n$
for $m\gg 0$ hence $\rho^{arb}(\gamma)=Id_n$ one concludes that $F=I_n$ or in other words that $\gamma$ lies in the
kernel of every representation $\pi_1(X)\to GL_n(A)$,  for an arbitrary $\C$-algebra $A$. In particular $\gamma\in H_n$. 
\end{pf}

\begin{coro}
 Assume $\pi_1(X,x)$ has a faithful representation in $GL_n(\C)$. Then $\widetilde{H_M^{\infty}}=\{ 1 \}$.
\end{coro}

\section{Rationality lemma}

 \subsection{Some pure Hodge substructures attached to an absolute closed  set $M$ and a fiber of $sh_M$}

 Let $f:Z\to X$ be a morphism  and $M\subset M_B(X,G)$ an absolute closed subset.

For ${\mathbb V}$ be an object
 of $\mathcal{T}_M$, we denote by $\text{tr}:\mathbb{V}\otimes \mathbb{V}^* \to \C$ the natural 
contraction.  Consider
the subspace $P_{\mathbb V}(Z/X)\subset H^1(Z,{\mathbb C})$
defined by:
\[
P_{\mathbb V}(Z/X) := \mathrm{Im}
\left[
f^{*}H^{1}(X,\mathbb{V})\otimes
  H^0(Z,{\mathbb V}^*) \stackrel{\cup}{\longrightarrow}
  H^{1}(Z,\mathbb{V}\otimes \mathbb{V}^{*})
  \stackrel{\text{tr}}{\longrightarrow} 
 H^1(Z,\mathbb{C})\right].
\]

In this formula, we denoted by $\V$ the local system on $Z$ defined as $f^*\V$. Obviously, no 
confusion can arise from this slight abuse of notation.

\begin{defi} We also define $P_M(Z/X), \overline{P}_M(Z/X) \subset H^1(Z,\C)$ as follows:

\begin{description}
\item[$P_M(Z/X) \subset H^1(Z,{\mathbb C})$:] the subspace of
  $H^1(Z,\mathbb{C})$ spanned by the $P_{\mathbb V}(Z/X)$, when
  ${\mathbb V}$ runs over all objects in $\mathcal{T}^{VHS}_M$.
\item[$\overline{P}_M(Z/X) \subset H^1(Z,{\mathbb C})$:] the subspace of
  $H^1(Z,\mathbb{C})$ spanned by the $P_{\mathbb V}(Z/X)$, when
  ${\mathbb V}$ runs over all objects in $\mathcal{T}_M$.
\end{description}
\end{defi}

$H^1(Z,\mathbb C)$ is defined over
$\mathbb Z$ since it is the complexification
of $H_{sing}^1(Z,\mathbb Z)$. This
Betti integral structure is the
one we will
tacitly use.

\begin{lem}
$P_M(Z/X)$  is a pure $\C$-Hodge
substructure of weight one of the $\C$-MHS underlying Deligne's
 MHS on $H^1(Z,\C)$.
 \end{lem}
 
 \begin{pf}
 Since each
  $\mathbb{V}$ is a $\mathbb{C}$-VHS of weight zero, and $X$ is
  smooth, it follows that $H^{1}(X,\mathbb{V})$ is a pure
  $\mathbb{C}$-Hodge structure of weight one. Also by \cite{Hod3} the mixed
  Hodge structures on the cohomology of varieties with coefficients in
  variations of Hodge structures are functorial and hence
  $P_{\mathbb{V}}(Z/X)$ is a $\mathbb{C}$-Hodge substructure of
  $H^{1}(Z,\mathbb{C})$.  Finally by strictness \cite{Hod3} the span
  $P_M(Z/X)$ of the $P_{\mathbb{V}}(Z/X)$'s will also be pure and of
  weight one.
 \end{pf}

\begin{lem}
If $G$ is defined over $\Q$ and that the absolutely closed subset $M\subset M_B(X,G)$ is defined over $\Q$, $\overline{P}_M(Z/X)$ is defined over $\mathbb Q$
\end{lem}

Assume now that $f(Z)$ is contained in a fiber of the reductive Shafarevich morphism for $M$ or that equivalently 
a finite \'etale cover of $Z$ lifts to a compact analytic subspace of $\widetilde{X_M}$. Then after a finite
etale cover we may assume that $f_*\pi_1(Z,z)\subset H_M$, ie that every object $\rho$ in $\mathcal{T}_M$
satisfies $\rho(\pi_1(Z,z))=\{ e \}$. 
The rationality lemma is the following statement: 

 \begin{theo}
 \label{rat} Assume $G$ is defined over $\Q$ and $M= M_B(X,G)$. Assume that $f_*\pi_1(Z,z)\subset H_M$.
 If $\pi_1(Z)\to \Gamma_M$ is trivial
 then $P_M(Z/X)=\overline{P}_M(Z/X)$. 
\end{theo}

\begin{coro}
 If $G$ and $M= M_B(X,G)$ are defined over $\Q$, then  $P_M(Z/X)$ is also defined over
 $\mathbb Q$.
\end{coro}

The rest of this section will be devoted to the proof of theorem \ref{rat}. 

 We will also assume $\dim M >0$
since the result is obvious for an absolute closed subset consisting of isolated points. The proof will be done in several steps which reduce the general statement to special situations. 

\begin{rem}
It seems likely that Theorem \ref{rat} holds true for arbitrary absolute closed subsets defined over $\Q$. One basically needs to adapt \cite{ES} to this situation.  
\end{rem}

\subsection{Reduction to the smooth case}

First we reduce to the case when $Z$ is
  smooth. We need the following lemma:

\begin{lem} $\overline{P}_M(Z/X)$
 is a pure weight one
substructure of Deligne's
 MHS on $H^1(Z)$.
\end{lem}

\begin{pf} Let ${\mathbb V}$ be an object
 of $\mathcal{T}_M$
By \cite[Theorem~4.1]{MTS} the space $H^{1}(X,\mathbb{V})$ carries a
pure twistor structure of weight one. Furthermore by
\cite[Theorem~5.2]{MTS} the space $H^{1}(Z,\mathbb{V})$ carries a
canonical mixed twistor structure and $f^{*}H^{1}(X,\mathbb{V})
\subset H^{1}(Z,\mathbb{V})$ is a twistor substructure. By
functoriality $P_{\mathbb{V}}(Z/X) \subset H^{1}(Z,\mathbb{C})$ will
be a pure weight one twistor substructure and  hence the span
$\overline{P}_M(Z/X) \sum_{\mathbb{V}}  P_{\mathbb{V}}(Z/X) \subset
H^{1}(Z,\mathbb{C})$ is a pure weight one 
twistor substructure of the mixed Hodge structure
$H^{1}(Z,\mathbb{C})$. However the  Dolbeault
realization of $\overline{P}_M(Z/X)$ is clearly preserved by 
$\mathbb{C}^{*}$ since by assumption $\mathbb{C}^*$ leaves $M^{Dol}$ invariant. Therefore $\overline{P}_M(Z/X)$ is a
sub Hodge structure. \ \hfill \end{pf}

In order to prove Theorem \ref{rat}, since $P_M(Z/X)\subset \overline{P}_M(Z/X)$ is
pure of weight one, it is
enough to prove that  $Gr_1^W P_M(Z/X)= Gr_1^W 
\overline{P}_M(Z/X)$. Hence,
  without a loss of generality, we can assume  that $Z$ is
smooth. 

\subsection{Reduction to a finite number of local systems}

 \begin{lem}\label{finiteness}
There is a  finite set $S$ of objects of
$\mathcal{T}^{VHS}_M$
 such that whenever a morphism $Z\to X$ has the property
 $\text{im}\left[\pi_1(Z)\to \Gamma_M\right] = 0$ it follows that 

\[
P_{M}(Z/X) = \sum_{\mathbb{V} \in S} P_{\mathbb{V}}(Z/X). 
\]

Similarly, there is a finite set $\overline{S}$ of objects of
$\mathcal{T}_{M}$, so that $\overline{P}_{M}(Z/X) = \sum_{\mathbb{V}
  \in \overline{S}} P_{\mathbb{V}}(Z/X)$. Furthermore the set
$\overline{S}$ can be chosen so that for any Higgs bundle $(E,\theta)$
corresponding to a $\V\in \overline{S}$ the $\C$-VHS associated to
$\lim_{t\to 0} (E,t.\theta)$ belongs to $S$.
\end{lem}

\begin{pf} Consider $(S_{\alpha})_{\alpha}$ a stratification of
 $Sh_M(X)$ by locally closed smooth algebraic
 subsets such that $s_{\alpha}:=sh_M|_{ (sh_M)^{-1}(S_{\alpha}) }:
 (sh_M)^{-1}(S_{\alpha}) \to S_{\alpha}$ is a topological fibration.
 Fix $p_{\alpha}\in S_{\alpha}$.
 Let $Z_{\alpha}=s_{\alpha}^{-1}(p_{\alpha})$, let
 $Z_{\alpha,o}$ be a connected component and let $Z'_{\alpha,o}\to Z_{\alpha,o}$ be the topological  covering space 
defined by
  $Z'_{\alpha,o}=
   \widetilde {Z^{univ}_{\alpha,o}}/
   \ker ( \pi_1(Z_{\alpha,o})\to \Gamma_M)$.
$Z'_{\alpha,o}\to Z_{\alpha,o}$.

Since $H^1(Z_{\alpha,o},\C)$ is finite dimensional,
it follows that a finite set $S$ exists with the required properties
for $Z= Z_{\alpha,o}$. Since the  cohomology classes coming from
$X$ are flat under the Gauss Manin connection,
 this statement holds true for all fibers
 of $s_{\alpha}$. Since every $f:Z\to X$ with the
 required properties factors through one of the $Z_{\alpha,o}$'s,
 the lemma follows.
\end{pf}

\subsection{Hodge theoretical argument}

From now on, we really need to assume that $M=M_B(X, G)$,
 and that $G$ is defined over $\Q$.

Let $A$ be a noetherian $\C$-algebra and $\rho_A: \pi_1(X,x) \to GL_N(A)$ be a representation. 
Let $\V_A$ be the
local system of free $A$-modules attached to $\rho_A$ and $\V_A^{\vee}=\mathrm{Hom}_A(\V_A, A)$ be the local system associated to $^t\rho^{-1}$.
We define:
$$P(A)= \mathrm{Im} \left[H^1(X,\V_A)\otimes_A H^0(Z, \V_A^{\vee}) \to H^1(Z,\C)\otimes_{\C} A \right]
$$
$P(A)$ is an $A$-submodule of the free $A$-module  $H^1(Z,\C)\otimes_{\C} A $.

Let $\sigma$ in $M^{**}$ be a non-isolated point. In subsection \ref{mhdef}, we recalled the construction
and basic properties
of $T_{\sigma} \subset R(\pi_1(X,x), G)$ a formal local subscheme which gives rise to a hull of the deformation functor of $\sigma$. In follows from \cite{GM2}, that this formal subscheme is actually the 
formal neigborhood of $\sigma$ in an analytic germ 
$T_{\sigma}^{an}\subset R(\pi_1(X,x), G)$. If  we decompose the reduced germ of $T_{\sigma}^{an}$ into the union
$T_{\sigma}^{an,red} = \cup_{i} T^{an,i}$ of its analytic irreductibile components, then we will denote by $T^i$
 the formal neighborhood of $\sigma$ in  $T^{an,i}$. The irreducible components of an analytic germ being 
in one to one correspondance with the irreducible components of the associated formal germ, it follows
that $T^{red}_{\sigma}= \cup_{i} T^i$ is 
still the irreducible decomposition of the reduced formal local subscheme underlying $T_{\sigma}$. 
Note that $T^i$ is an integral formal
  subscheme of $T_{\sigma}$ and so its 
ideal $\mathfrak{P}^i$ is a minimal prime  of
$\widehat{\mathcal{O}}_{T_{\sigma}}$.

\begin{lem} \label{mhprime} The weight and Hodge filtrations on $\widehat{\mathcal{O}}_{T_{\sigma}}$ induce on
 $\mathfrak{P}^i\subset \widehat{\mathcal{O}}_{T_{\sigma}}$ a sub-MHS structure. 
\end{lem}
\begin{pf}
First observe that the minimal associated primes 
of the graded ring $Gr^{\m^{\bullet}} \widehat{\mathcal{O}}_{T_{\sigma}}$
are graded ideals and also split subMHS of the split MHS on $Gr^{\m^{\bullet}} \widehat{\mathcal{O}}_{T_{\sigma}}$
since the $Res_{\C|\R}\C^*$-action defining the Hodge decomposition is compatible with the ring structure. 

 By construction, 
there is a ring isomorphism $\widehat{\mathcal{O}}_{T_{\sigma}} \to Gr^{\m^{\bullet}} \widehat{\mathcal{O}}_{T_{\sigma}}$. This ring isomorphism takes minimal associated primes to minimal 
associated primes.  Hence, $Gr^{\m^{\bullet}}\mathfrak{P}^i$ is a sub Hodge Structure of $Gr^{\m^{\bullet}} \widehat{\mathcal{O}}_{T_{\sigma}}$.

There is no canonical choice for this isomorphism but it can be chosen in such a way that it respects the weight and Hodge filtrations - but not the three filtrations. This implies that the 
trace of  the Hodge  filtration of $Gr^{\m^{\bullet}} \widehat{\mathcal{O}}_{T_{\sigma}}$ on  
$Gr^{\m^{\bullet}}\mathfrak{P}^i$ is the filtration induced by the trace of the Hodge filtration of $\widehat{\mathcal{O}}_{T_{\sigma}}$ on $\mathfrak{P}^i$. The anti Hodge filtration satisfies a similar statement.

These two facts imply that $\mathfrak{P}^i\subset \widehat{\mathcal{O}}_{T_{\sigma}}$ is a  sub-MHS structure.
\end{pf}

 Hence the complete local algebra $\widehat O_{T^i}$ carries a $\C$-MHS and $\rho_ {\widehat O_{T^i}}:\pi_1(X,x)\to G(\widehat O_{T^i})$ is 
the monodromy of the local system
$ \D(\V_{\sigma})\otimes_{\widehat O_{T_{\sigma}}} \widehat O_{T^i}$. Thanks to lemma \ref{mhslice} and lemma
\ref{mhprime},  this local system underlies a $\C$-VMHS whose weight filtration corresponds to the powers of the maximal ideal in $\widehat O_{T^i}$.

By, construction the tautological
  representation $\rho_{\mathcal{O}_{T^{an,i}}}: \pi_1(X,x) \to 
G(\mathcal{O}_{T^{an,i}})$ is a holomorphic family of representations
parametrized by a reduced germ of complex space.

If there
is a proper closed analytic subset $Z^i\subset T^{an,i}$ such that
$\forall p\in T^{an,i}-Z^i$ the representation
$\rho_{\mathcal{O}_{T^{an,i}}}(p)$ is a reductive representation, then
the inclusion $f_*\pi_1(Z,z)\subset H_M$ implies that the restriction
of $\rho_{\mathcal{O}_{T^{an,i}}}(p)$ to $\pi_1(Z,z)$ is trivial for
$p\not\in Z^i$.
Hence the restriction of $\rho_{O_{T^{an,i}}}$ and $\rho_{\widehat O_{T^i}}$ to $\pi_1(Z,z)$ are trivial as well. 

If not, then for each irreducible component $M'\subset M$ containing $\sigma$, take a component $R'$ via $\sigma$ of the preimage $\pi^{-1}(M')\in R(\pi_1(X,x),G)$ which dominates $M'$.
Let $(R')^{red} \subset R'$ be its maximal reduced subscheme. Consider the semisimplification of  the
representation attached to the generic point of the subscheme $(R')^{red} \subset R(\pi_1(X,x),G)$. It is conjugate to a Zariski dense representation with values in with values
in some $G'\subset G$, where $G'$ is reductive over $\overline{\Q}$.
But
$\mathrm{Im}(M_B(X,G')\to M_B(X,G))$ is a closed acqc set and so $M'
\subset \mathrm{Im}(M_B(X,G')\to M_B(X,G))$. So without a loss of
generality we may replace $G$ by $G'$ and also replace
$T_{\sigma}^{an}$ by an analytic Goldman-Millson slice through
$\sigma$ in $R(\pi_1(X,x), G')$.
 With this new definition, the restriction of $\rho_{\widehat O_{T^i}}$ to $\pi_1(Z,z)$ is trivial too
and the corresponding local system on $Z$
is the constant local system $\V_{\sigma} \otimes_{\C} \widehat O_{T^i}$. 

In particular, we have a canonical isomorphism of VMHS $\V_{\widehat O_T^i/ \m^k}^{\vee}|_Z\simeq \V_{\sigma}^{\vee}|_Z\otimes _{\C} \widehat O_T^i/ \m^k$. It now follows that, for all $k\in \N$:

$$P(\widehat O_{T^i}/ \m^k)= \mathrm{Im} (H^1(X,\V_{\widehat O_T^i/ \m^k})\otimes_{\C} H^0(Z, \V_{\sigma}^{\vee})\buildrel{H_{\widehat O_T^i/ \m^k}}\over{\to} H^1(Z,\C)\otimes_{\C} \widehat O_{T^i} / \m^k).
$$

$H_k:=H_{\widehat O_{T^i}/ \m^k}$ preserves the natural Mixed Hodge structures.

\begin{prop} \label{mharg} $P_k:= P(\widehat O_{T^i}/ \m^k)\subset P_1\otimes \widehat O_{T^i}/ \m^k \subset P_M(Z/X)\otimes \widehat O_{T^i}/ \m^k  $.
\end{prop}

\begin{pf} If $k=1$ this is trivial: 
by construction, $P_1\subset P_M(Z/X)$. We now argue by induction and assume that the result holds for $k'<k$. 

The representation $\rho_k=\rho_{\widehat O_T^i/ \m^k}$  underlies a variation of complex mixed Hodge structure $\M_k$
on $X$. The weight filtration is given by the powers of $\m$. 
Since $\rho_k$ is trivial on $\pi_1(Z)$ then its restriction to
$Z$ is the trivial VMHS $\mathbb H \otimes_ {\mathbb C} \widehat O_T^i/ \m^k$
where $\mathbb H$ is some Hodge structure of weight zero (with a possibly non trivial Hodge vector) and, on $\widehat O_{T^i}/ \m^k$, the weight filtration is described by the powers of $\m$.

\[
P_k= \mathrm{Im} \left[\xymatrix@1@C+1pc{H^1(X,\M_k\otimes_{\C} \H)
    \ar[r]^-{H_{k}} & 
H^1(Z,\C)\otimes_{\C} \widehat{\mathcal{O}}_{T^i} / \m^k}\right].
\]

The weights of $\M_k$ are $0, \ldots, -k+1$. 
Consider  the following diagram of $MHS$, in which the rows are exact:

\[
\xymatrix{
H^1(X, W_{-k+1} \M_k\otimes_{\C} \H ) \ar[r] \ar[d] & H^1(X, \mathbb
  \M_k\otimes_{\C} \H) \ar[r] \ar[d] & H^1(X, \M_{k-1}\otimes_{\C} \H)
  \ar[d] \\
H^1(Z) \otimes \m^{k-1} /\m^k \ar[r]  &
H^1(Z) \otimes 
  \widehat{\mathcal{O}}_{T^i} / \m^k \ar[r] & H^1(Z)\otimes
  \widehat{\mathcal{O}}_{T^i} / 
  \m^{k-1}
}
\]

Remember we assume $Z$ to be smooth.
The weights of the MHS in the first row are $2-m$,  in the second
$2-m, ..., 1$, in the third one $3-m, ..., 1$. Hence the second line
is just the canonical exact sequence 

\[
\xymatrix@1{
W_{2-k}\left[
H^1(Z) \otimes \widehat{\mathcal{O}}_{T^i} /
  \m^k
\right] 
  \ar[r] & H^1(Z)
  \otimes \widehat{\mathcal{O}}_{T^i} / \m^k \ar[r] & 
 Gr^W_{3-k}\left[ H^1(Z) \otimes \widehat{\mathcal{O}}_{T^i} /
  \m^k\right].
}
\]

The main observation is now that, by strictness, we have:

\[
W_{2-k} P_k= Im\left[ 
\xymatrix@1{
H^1\left(X, W_{-k+1}\left(\M_k\otimes_{\C} \H\right)\right) \ar[r] & 
H^1(Z)
\otimes \widehat{\mathcal{O}}_{T^i} / \m^k}\right].
\]

From this it follows that $ W_{2-k} P_k\in P_1\otimes \m^{k-1}  / \m^k$. 
By induction, $P_ {k-1}\subset H^1(Z)\otimes \widehat O_{T^i} / \m^{k-1}\subset P_1 \otimes \widehat O_{T^i} / \m^{k-1}$.
But $P_k$ (respectively $P_{k-1}$) is the image of the map in the third column (resp. the second). It follows that $P_k\subset P_1\otimes \widehat{O_{T^i}} / \m^k$.

\end{pf}

\subsection{Proof of theorem \ref{rat} if $M=M_B(X,G)$}

It follows from proposition \ref{mharg} that:

\[
P(\mathcal{O}_{T^{an,i}}) \subset P_M(Z/X)\otimes \mathcal{O}_{T^{an,i}}.
\]

It follows that for all $p$ in the complex analytic germ $T^{an,i}$ we have:
$$P_{\mathbb V_{\rho(p)}}(Z/X)\subset P_M(Z/X).$$

Since there is a
  complex analytic neighborhood $U$ of $\sigma$ in $M$ such that every
  point of $U$ has a (semisimple) representative in $T^{an,i}$, it
  follows that for every $\V\in U$ we have $P_{\mathbb V}(Z/X)\subset
  P_M(Z/X)$.

Now let $\overline{S}$ be
  the finite set from Lemma~\ref{finiteness}. Suppose $\V \in
  \overline{S}$ with an associated Higgs bundle $(E,\theta)$, and let
  $\left(\V_t\right)_{t\in\C^*}$ be the local systems corresponding to
  the Higgs bundle $(E,t\theta)$. For a small enough $t$ we have:
$$
P_{\V_t}(Z/X)\subset P_M(Z/X)
$$

Fix $t$ small enough and non zero. It follows that
 $\dim(\sum_{\V\textit{}\in \bar S} P_{\V_t}(Z/X))\le \dim P_M(Z/X)$.

Consider: 

\[
P^{Dol}_{\mathbb V_t}= \mathrm{Im}\left[
\xymatrix@1@C+1pc{H^1_{Dol}(X, \V_t)\otimes
H_{Dol}^0(Z,{\mathbb V}_t^* ) \ar[r]^-{{\mathrm{Id}}\otimes
{\mathrm{tr}}} &  H_{Dol}^1(Z)
}
\right].
\]

Using Simpson's Dolbeault isomorphism we have:

\[
\dim\left(\sum_{\V\in \bar S} P^{Dol}_{\V_t}(Z/X)\right)\le \dim
P_M(Z/X).
\]

Recall there is a natural isomorphism $s(t):H^{\hdot}_{Dol}(-,\V) \to H^{\hdot}_{Dol}(-,\V_t)$. Let $(E,\theta)$ be a polystable Higgs bundle representing $\V$. Then $H_{Dol}^{\bullet}(X, \V):= \mathbb{H}^{\bullet} (X, (E\otimes \Omega^{\bullet}_X, \theta))$. We can construct an quasi-isomorphism $(E\otimes \Omega^{\bullet}_X, \theta)\to (E\otimes \Omega^{\bullet}_X, t.\theta)$ by the formula:

\[
\xymatrix@1@C+1pc{
E \ar[r]^{\theta} \ar[d] & E\otimes \Omega^1_X \ar[r]^{\theta\wedge} \ar[d]^{t} & E\otimes \Omega^2_X \ar[r]^{\theta\wedge} \ar[d]^{t^2} &\ldots \\
E \ar[r]^{t\theta}  & E\otimes \Omega^1_X \ar[r]^{t\theta\wedge}  & E\otimes \Omega^1_X \ar[r]^{t\theta\wedge}  &\ldots
}
\]

Since $(E,t\theta)$ is the polystable Higgs bundle representing $\V_t$ this quasi isomorphism defines indeed an isomorphism $s(t):H^{\hdot}_{Dol}(-,\V) \to H^{\hdot}_{Dol}(-,\V_t)$. 

In case $(E,\theta)$ is kept fixed by $\C^*$
which means that there is an isomorphism $\psi(t): (E,\theta)\to (E, t\theta)$, $a(t)=\psi(t)^{-1}\circ s(t)$
is an automorphism of $H^1_{Dol}(X, \V)$ which comes from an action of $\C^*$. Here $(E,\theta)$ is not kept fixed by $\C^*$ but its restriction to $Z$ is. This gives a diagram:

\[
\xymatrix{
H^1_{Dol}(X, \V) \ar[r] \ar[d]^{s(t)} & H^1_{Dol}(Z, \V)  \ar[d]^{a(t)} \\
H^1_{Dol}(X, \V_t) \ar[r]  & H^1_{Dol}(Z, \V) \
}
\]

By functoriality and the definition of $P^{Dol}_{\V_t}(Z/X)$, we get a commutative diagram:

\[
\xymatrix{
P^{Dol}_{\V}(Z/X) \ar[r] \ar[d]^{\tilde s(t)} & H^1_{Dol}(Z)  \ar[d]^{a(t)} \\
P^{Dol}_{\V_t}(Z/X) \ar[r]  & H^1_{Dol}(Z ) \
}
\]

and $\tilde s(t)$ is an isomorphism. 

Hence $\dim(\sum_{\V\in \bar S} P^{Dol}_{\V}(Z/X))\le \dim P_M(Z/X)$. Since, by Simpson's Dolbeault isomorphism,  the l.h.s is $\dim \overline{P}_M(Z/X)$ the
theorem is proved.

 \section{Construction of the Shafarevich morphism}

\subsection{Preliminary considerations}

\subsubsection{Pure weight one rational subspaces of
  $H^1 ( Z )$  }

Let $Z$ be a complex projective variety. 

The possibly non zero Hodge numbers
 of Deligne's Mixed Hodge
structure
 \cite{Hod3} on the first cohomology group $H^1(Z,{\mathbb Z})$
 of the connected projective
 variety $Z$
are $h^{0,0},h^{0,1},h^{1,0}$.

In particular, we have  an extension of
$\mathbb Q$-MHS of a pure weight one
HS by a pure weight zero HS:
\begin{equation}
\label{eqn0} 0 \to W_0 (H^1(Z,{\mathbb Q}))
\to H^1(Z,\mathbb Q)
\to Gr_1^W (H^1(Z,{\mathbb Q}))\to 0.
\end{equation}

Let $Z^{sn}\to Z$ be the seminormalisation
of $Z$
(see \cite{Kolrat} Chap. I Definition 7.2.1,
p. 84 and the original references therein)
 $H^1(Z)\to H^1(Z^{sn})$ is
 an isomorphism of MHS since
 $Z^{sn}(\mathbb C)\to Z(\mathbb C)$
is an homeomorphism \cite{Kolrat} I.(7.2.1.1).

Let $\mathcal A$ be a pure weight
 one $\mathbb Q$-HS. There is an abelian variety
(well defined up to isogeny) such that
 $H^1(A, {\mathbb Q}) ={\mathcal A}$.

\begin{lem} \label{jacobi} Let
 $\widetilde{\phi}:
 \mathcal A \to H^1(Z,{\mathbb Q})$ be
a morphism of MHS. Then there exists
 a rational number
$d\not=0$ and a
morphism $\psi:Z^{sn} \to A$ such that
$d\widetilde{ \phi}= H^1(\psi)$.
\end{lem}

We may as well assume $Z$ is
seminormal. Assume moreover that $Z$ is a curve.
Consider more generally
 $\phi: {\mathcal A} \to Gr_1^W H^1(Z)$ 
a  morphism  $\mathbb Q$-HS 
of pure weight one.

 Pulling back
the extension \eqref{eqn0} 
by the morphism
$\phi$ defines
 a  extension of
$\mathbb Q$-MHS 

\begin{equation}\label{ext}0 \to W_0 (H^1(Z,{\mathbb Q}))
\to {\mathcal A}'
\to {\mathcal A} \to 0.
\end{equation}

\begin{pf}
Let $\nu:Z^{\nu} \to Z$ be
the normalisation of $Z$.
 Thanks to \cite{Hod3} lemme 10.3.1.
 the extension \eqref{eqn0}
 is isomorphic to
 $$0 \to W_0 \to H^1(Z)
\buildrel{\nu^*}\over{\longrightarrow}
 H^1(Z^{\nu})\to 0.$$

Let $\gamma:[0,1]\to Z$ be a loop
which is
based at a singular point, meets the singular
locus of $Z$ at finitely many points and
is smooth outside these points. 
The preimage of $\gamma$
in $Z^{\nu}$
$\gamma$ is a finite
union $\gamma_1,\ldots,\gamma_n$
of paths possibly lying in several
connected components of $Z^{\nu}$.
This defines a linear form
$\int_{\gamma}:\omega \mapsto
 \sum_i \int_{\gamma_i}\omega$
on $H^0(Z^{ \nu},\Omega^1)$ and, upon composition with
 $\phi$, a linear form $\phi^*\int_{\gamma}$ on
${\mathcal A}^{1,0} $.

It follows from
\cite{Carext} theorem (1.13)
 -see also the enlightening
example (1.17)
- that \eqref{ext}
is split if and only if
 for every $\gamma$ as above
 $\phi^*\int_{\gamma}$
is a rational multiple of a period of
$\mathcal A$, i.e. lies in the
 image of $H_1(A,{\mathbb Q})$.

The datum $\widetilde{\phi}$
 gives actually such a
 splitting and the Abel-Jacobi construction 
gives a continuous mapping $Z\to A$ with the
required property which is holomorphic when 
pulled back to $Z^{\nu}$. Since $Z$ is
 seminormal this continuous mapping actually
underlies a morphism.

The general case readily follows
from the curve case.
 Assume
first $Z$ is irreducible. Let $\lambda: Z^{\nu}\to Z$ be the normalisation of $Z$. Then we can construct
a morphism
$\psi^{\nu}: Z^{\nu} \to A$ and an integer $d$ such that $d H^1(\lambda)\circ \tilde {\phi } = H^1(\psi)$ . This morphism  is locally constant
 on the fibers of $\lambda:Z^{\nu}\to Z$.
 On the other hand we can always find a connected curve $C$
 passing through each connected component of a given positive dimensionnal fiber $F$
 of $\lambda$. Consider $C^{sn}\to C$ the seminormalization of $C$. This is a homeomorphism
which identifies $H^1(C)$ and $H^1(C^{sn})$ with their respective Mixed Hodge structures. 
 The morphism $\psi^{\nu}|_C:C^{sn}\to A$
  is isogenous to the one predicted by
 lemma \ref{jacobi} applied to $C$ and the resulting $\tilde \phi_C:
 \mathcal{A}\to H^1(C)$. Let $C'$ be the image of $C$ in $Z$.
 Since $\tilde \phi_C$ factors through $H^1(C')$ it follows
 that $\psi^{\nu}|_C$ is constant on the finite fiber of $C^{sn} \to C'$.
Hence $\psi^{\nu}$ assume the same value on all connected components of $F$.
Hence it descends to morphism $\psi: Z\to A$ since $Z$ is seminormal. 

 In general,
$Z$ has $m$ irreducible components,
there are $m-1$ constants of integration to
take care of and a connected curve in $Z$
 meeting every connected component
 of the smooth locus do the bookkeeping.
 \end{pf}

 \subsubsection{Period mappings for $\C$-VMHS}

 $\R$-MHS have period domains and  $\R$-VMHS period mappings generalizing 
 those constructed by Griffiths for $\R$-VHS, \cite{U},  see also \cite{Carext} \footnote{Actually, $\C$-VMHS have also period
 mappings of their own but since this would  not give additional information, we will stick to the usual conventions
 used in the litterature}. 

Recall that $X$ is a complex projective manifold and
let   $(X, \mathbb{V},\mathcal{F}^{\bullet}, S)$ be a $\R$-VHS of weight zero  and let $M$ be the real Zariski closure of its monodromy group computed at some basepoint $x\subset X$. 
Let $U \subset M$ be the
  isotropy group of the Hodge filtration on $\mathbb{V}_{x}$.
Then the period
  domain of $\mathbb{V}$ is the complex manifold $D(\mathbb V) :=
  M/U$. It is endowed with a certain horizontal distribution which can
  be described in terms of the Hodge structure on the Lie algebra
  $\mathfrak{m}$ of $M$.   It is actually the actually the 
period domain attached to the Hodge semisimple group $M^{ad}$. Furthermore $D(\mathbb V)$ is a moduli space of Hodge structures on 
$M$, see \cite{GS} for more details. 
 
 Let  $(X, \mathbb{V},\mathbb{W}_{\bullet},\mathcal{F}^{\bullet}, (S_k)_{k\in\Z})$ be a  $\R$-VMHS.
Again we have a
   period domain $MD(\mathbb{V})$ for this variation and a holomorpic
   fibration of period domains $\psi: MD(\mathbb{V})\to \prod_k
   D(Gr_{k}^{\mathbb{W}} \mathbb{V})$ which is compatible to the
   horizontal distributions.  

The domain $MD(\mathbb{V})$ is a
   homogenous space of the form $ H/U'$, where $H$ is the subgroup of
   $W_0 GL( \mathbb{V}_x)$ mapping to $\prod_k M(S_k)$ under the
   natural surjection $W_0 GL( \mathbb{V}_x)\to GL( Gr_0^ {\mathbb W}
   \mathbb{V}_x)$.

 Accordingly there is an  equivariant holomorphic horizontal period mapping $\phi_{\mathbb{V}} : \widetilde{X^{univ} }\to MD( \mathbb{V})$ with a commutative diagram:
 
  $$
  \begin{array}{ccc}
  \widetilde{X^{univ}} & \buildrel{\phi_{\mathbb{V}}}\over{\longrightarrow} & MD(\mathbb{V}) \\
  &\searrow & \downarrow  \psi\\
  & &  \prod_k D(Gr_{k}^{\mathbb{W}} \mathbb{V}).

  \end{array}
  $$

 Let  $\mathbb M$ be a $\R$-VMHS of weights -1 , 0
  and $MD$ be its period domain.
 Let $D$ be product of the the period domains corresponding
  to the graded parts of $\mathbb{M}$. The map $MD\to D$ is then an affine bundle. 
  
The following lemma can be extracted from
  \cite[p.~200]{Carext}.
  \begin{lem}

  $MD\to D$ is a holomorphic vector bundle.

  The fiber $V (H_{-1},H_0)$ of $MD\to D$ at $(H_{-1},H_0)$
  is canonically isomorphic to $\mathrm{Hom} (H_0, H_{-1})_{\C}/ F^0$ where 
  $\mathrm{Hom} (H_0, H_{-1})$ is endowed of its natural Hodge structure of weight $-1$.
  \end{lem}

  Consider $f:Z\to X$
  a morphism such that $f^*Gr_i \mathbb{M}$ is a VHS with trivial monodromy. Let $P_{\mathbb{M}}= \sum_i P_{f^*Gr_i \mathbb{M}} (Z/X)\subset H^1(Z)$, then we have the following lemma:

  \begin{lem} There is a commutative diagram
  $$
  \begin{array}{ccc}
  \widetilde{Z^{univ}} & \to & (P_{\mathbb{M}}^{1,0})^* \\
  &\searrow & \downarrow  g_{\mathbb{M}}\\
  & & V (H_{-1},H_0)

  \end{array}
  $$
where $g_{\mathbb{M}}$ is linear and injective and, when $Z$ is smooth, the horizontal map is given by integration of closed holomorphic forms.
  \end{lem}

\begin{pf}
The proof is straightforward and is left to
  the reader. The case when  $\V=\C$ is standard and the general case
  follows by the same reasoning.
\end{pf}

\subsection{Proof of Theorem \ref{main}}

  \subsubsection{Notations}

In what follows, $M= M_B(X, G)$ where $G$ is a reductive group defined over $\Q$. 

  \begin{lem} 
  \label{choice}
There exists an object
    $\mathbb{M}_{1}$ of $\mathcal{T}_M^{MVHS}(1)$ such that for every
    $f:Z\to X$ for which $\pi_1(Z)\to \Gamma_M$ is trivial, we have
    that $P_{\mathbb{M}_1}=P_M(Z/X)$ and that $g_{\mathbb{M}_1}$ is
    injective
  \end{lem}
\begin{pf}
Take 
\[
\M_1 :=
  \sum_{\sigma \in S} \left(\D_{1}(\V_{\sigma})+
  \overline{\D_{1}(\V_{\sigma})}\right),
\] 
where $S \subset \mathcal{T}_M^{VHS}$ is the finite set constructed in
lemma~\ref{finiteness}, and $\D_1(\V_{\sigma})$ is the
$\mathbb{C}$-VMHS from Definition~\ref{defi:Dk}. 
\end{pf}

Let $\widetilde{X^k_M}$ be the covering space of $X$ defined as
$\widetilde{X ^{univ}}/ \widetilde{H^k_M}$.
This covering is Galois with Galois group
$\Gamma_M^k$.

Consider the local systems
  that belong to the finite set $S$ in $\mathcal{T}^{VHS}_M$ from lemma
\ref{finiteness}. Without loss of generality we may assume that they
underly real VHS of weight zero. Every $\rho$ in $S$
underlies a Zariski dense representation
$\pi_1(X) \to G_{\rho}$ where $G_{\rho}$
is a real Lie group of Hodge type. Let
$\rho_S: \pi_1(X)\to G_S= \prod_{\rho \in S} G_{\rho}$ be the direct sum
representation.

\subsubsection{Construction of the Shafarevich morphism in case $k=1$}

In this paragraph, we assume that $k=1$.
Choose a
finite dimensionnal real representation as in lemma \ref{choice}
of $\mathcal{G}_S^1(\mathbb R)$ such that the associated
local system $\mathbb {W}(1)$ underlies a graded polarizable real
variation of mixed Hodge structure with the finite weight filtration:
$0=\mathbb{W}_{-2}\subset \mathbb{W}_{-1} \subset
\mathbb{W}_0=\mathbb{W}(1)$.

 Associated with $\mathbb{W}(1)$,
 we have a holomorphic Griffiths' transversal
  period mapping
 $q^1_S:\widetilde {X^{univ}}\to \mathcal{D}^1_S$
 where $\mathcal{D}^1_S$ is is the corresponding period domain for MHS.
 The period domain $\mathcal{D}^1_S$ has a holomorphic
 fibration $\pi:\mathcal{D}^1_S \to \mathcal{D}_S $
 which makes it an affine fibration over the period domain
 $\mathcal{D}_S $.
The composition $\pi
 \circ q^1_S$ is the period mapping for the associated graded object of
 $\mathcal{T}^{VHS}_M$.

The map $q^1_S$ factors through a holomorphic horizontal map
$Q^1_M\widetilde{X^1_M }\to \mathcal{D}^1_S $.

 Consider the  holomorphic map
  $q_S:\widetilde{X^1_M}
   \buildrel{Q^1_M\times Sh_M}
   \over{\longrightarrow} \mathcal{D}^1_S \times Sh_M(X)$.

   \begin{lem}\label{ml}
    Every connected component of a fiber of $q$
   is compact.
   \end{lem}
   
\begin{pf}
Such a component $\Phi$ is contained in the lift of some fiber $Z$
of $X\to Sh_M(X)$. Replacing $Z$ by an etale cover,  we may
assume $\rho(\pi_1(Z))=\{ e\}$ whenever
 the conjugacy class of the reductive representation 
 $\rho$ is in $M$.
 
 Hence $\Phi$ is a connected component of a fiber
 of the map $q'$ defined as $q_S$ restricted to $\widetilde{Z^1_M}=
 \widetilde{Z^{univ}}/ \ker(\pi_1(Z)\to \widetilde{\Gamma^1_M})$.
 
 Now $\pi \circ q'$ is the constant map 
 and $\Phi$ is a connected component
 of a fiber of an holomorphic map
 $\psi : \widetilde{Z^1_M} \to V$ where $V$ is a
 complex vector space which is a fiber of
 $\pi$.
 
 Apply lemma \ref{jacobi} to $X=Z$ and $\mathcal{A}= P_M(Z/X)$.
The rationality hypothesis is fulfilled thanks to Theorem \ref{rat}. 
 We find a map to an abelian variety $Z^{sn} \to A$ and 
 using $P^{1,0}_M(Z/X)^*\to A$
 the universal covering space of $A$
 a proper holomorphic map 
 $\psi':\widetilde{Z^1_M}^{sn} \to P^{1,0}_M(Z/X)^*$.

 Our claim follows from the fact that we have a
 commutatative diagram:
 
 $$
 \begin{array}{ccc}
\widetilde{Z^1_M}^{sn}
&\buildrel{\psi}\over{\longrightarrow}&P^{1,0}_M(Z/X)^*\\
s\downarrow & & i\downarrow \\
\widetilde{Z^1_M} &\buildrel{\psi}\over{\longrightarrow} & V
 \end{array}
 $$

 Where $s$ is the seminormalisation and $i$ an injective linear map.

\end{pf}

Next, recall the following classical result:

\begin{lem} \label{locpropre}
(\cite{Car}, vol 2, pp.
797-811) Let $X, S$ be two
complex spaces
 and $f:X\to S$ a morphism. Assume a  connected
component  $F$ of a fiber of $f$  is compact. Then,
$F$
has a neighbourhood $V$ such that $g(V)$ is a local analytic
subvariety of $S$ and $V\to g(V)$ is proper.

Assume furthermore any  connected
component   of a fiber of $f$
  is compact and $X$ and $S$ are normal.
Then, the set $\bar S$ of connected
components of a fiber of $f$
can be endowed with a structure of normal complex space such that
the quotient mapping
$e: X\to \bar S$ is holomorphic, proper, with connected fibers.
\end{lem}

Using this lemma, we construct a surjective proper holomorphic mapping
with connected fibers to a normal complex
space $r^1_M:\widetilde{X^1_M}\to
\widetilde{S^1_M}(X)$ such that its fibers are precisely
 the connected components of the fibers of $q$.
 Since $q$ is $\Gamma^1_M$-equivariant it folllows
 that $r^1_M$ is $\Gamma^1_M$-equivariant too.
Note that $\Gamma^1_M$ acts on $\widetilde{S}^1_M(X)$ in a proper
discontinuous fashion and hence has at most finite stabilizers.

\begin{lem} The fibers of $r^1_M$ are precisely the maximal connected
analytic subvarieties of $\widetilde{X^1_M}$.
\end{lem}
\begin{pf}
 It is enough to show that whenever $Z$ is a connected
 compact analytic subvariety of $\widetilde{X^1_M}$,
  $r^1_M$ is constant. Fix such a $Z$.

 The map $f:Z\to X$ has the property that the group 
 homomorphism $\pi_1(Z)\to \Gamma^1_M$ induced by
 $\pi_1(f)$ has finite image.
Let $Z'$ be a connected \'etale cover of $Z$
such that $\pi_1(Z')\to \Gamma^1_M$ is trivial.
Abusing notation, let $f:Z'\to X$ be the resulting map.
Then, for every representation $\rho$ in $M$,
$f^*\rho$ is trivial and for every object $\mathbb V$ of
 $\mathcal{T}^{VHS}_M$, the restriction map
 $H^1(X,{\mathbb V} )
 \to 
 H^1(Z',{\mathbb V})$ is zero.
 This implies, through the proof of lemma \ref{ml}
 that $q$ is constant on $Z'$ and thus $r^1_M$.
 \end{pf}

{\bfseries Remark.}  In fact it can be shown that  $\widetilde{S^1_M}X)/ \Gamma^1_M$
 is a normal algebraic variety. This
follows from recent work of G. Pearlstein but is not used in the the
main theorem and so we will not discuss it here.

 \subsubsection{Stein property in the case $k=1$}
 
 \begin{prop}
  $\widetilde{X^1_M}$ is holomorphically convex
 and $r^1_M$ is its Cartan-Remmert factorisation.
 \end{prop}
 \begin{pf}
 Consider the natural period mapping  $\widetilde{S_M}(X)\to \mathcal{D}_S$
 and the affine bundle
  $V_S(X)=\widetilde{S_M}(X)\times_{ \mathcal{D}_S}
  \mathcal{D}^1_S \to \widetilde{S_M}(X) $.
  The previous consideration imply that
   $\widetilde{S^1_M}(X)\to V_S(X)$ is proper
    and finite to one, hence finite.

    Being an affine bundle over a Stein space $V_S(X)$ is Stein, hence
 $\widetilde{S^1_M}(X)$ is Stein.
\end{pf}
 
\subsubsection{General case}

 \begin{theo}
 Let $\widetilde{X}=\widetilde{X^{univ}}/\Gamma$ be a 
 Galois  covering space
 of $X$ with $\widetilde{H^{\infty}_M}\subset \Gamma\subset \widetilde{H^1_M}$.
 Then $\widetilde{X}$ is holomorphically convex.
 \end{theo}

 \begin{pf}
 Consider the map $q:\widetilde X\to \widetilde{S^1_M}(X)$.
 
 We claim that
    every connected component $\Phi$ of a fiber of $q$
   is compact. Indeed $\Phi$ has to be a connected lift
   of a projective variety $Z\subset X$ which is mapped to 
   a point in $S^1_M(X)$. Replacing $Z$ by an \'etale cover,
   we may assume $\pi_1(Z)\to \Gamma^1_M$ is trivial
   hence $Im(\pi_1(Z)\to\pi_1(X))\subset \Gamma$
    by Proposition \ref{strict}. This implies that $\Phi$ is compact.
    
    In particular we may construct its Stein factorization
    $\widetilde {X} \to \widetilde S$ and it follows from the 
    previous argument
    that $p:\widetilde{S} \to \widetilde{S^1_M}(X)$ has the following
    property: 
    
    \begin{lem}
    Every point $x\in\widetilde{S^1_M}(X)$ has a neighbourhood $U$
    such that $p^{-1}(U)$ is the disjoint union of open sets
    $V$ and $p|_V$ is a quotient map by a finite group $G$.
    \end{lem}
    
    This certainly implies that $\widetilde{S}$ is Stein.
    
    In fact the finite group in question 
    is $\ker (\pi_1((r^1_M)^{-1}(x))\to \Gamma/\Gamma^1_M)$
    and
    injects into the real points
    of
    a prounipotent proalgebraic group. It is
    thus a trivial group
     hence $\widetilde{S}\to\widetilde{S^1_M}(X)$
    is a topological covering map.
    \end{pf}

\vskip .3cm

Philippe Eyssidieux. 

Institut Fourier. Universit\'e Joseph Fourier. Grenoble (France). 

eyssi@fourier.ujf-grenoble.fr

\vskip .2 cm

Ludmil Katzarkov. 

Department of Mathematics. University of Miami. Miami (Florida, USA). 

Fakult\"at f\"ur Mathematik. Universit\"at Wien. Wien (\"Osterreich). 

lkatzark@math.uci.edu

l.katzarkov@math.miami.edu

\vskip .2 cm

Tony Pantev.

Department of Mathematics. University of Pennsylvania. Philadelphia (Pennsylvania, USA).

tpantev@math.upenn.edu

\vskip .2 cm

Mohan Ramachandran. 

Department of Mathematics. The State University of New York at Buffalo. Buffallo (NY, USA). 

ramac-m@math.buffalo.edu
\end{document}